\documentclass[a4paper, fontsize=12pt]{article}

\renewenvironment{abstract}
               {\list{}{\rightmargin\leftmargin}%
                \item[\hspace*{1cm}\small\textbf{Abstract ---}]\relax}
               {\endlist}

\usepackage[top=0.8in, bottom=0.8in, left=0.8in, right=0.8in]{geometry}
\usepackage[fleqn]{amsmath}
\usepackage{amssymb}
\usepackage{enumerate}
\usepackage{amsthm,bm}
\usepackage{appendix}
\usepackage[russian,english]{babel}
\usepackage{subfigure}
\usepackage{float}
\usepackage[pdftex]{graphicx}
\usepackage{framed}
\usepackage{hyperref}
\usepackage{sidecap}
\usepackage{mathrsfs}
\usepackage{xcolor}

\sidecaptionvpos{figure}{t}
\newtheorem{Theorem}{Theorem}[section]

\newtheorem{Corollary}[Theorem]{Corollary}

\newtheorem{Definition}[Theorem]{Definition}

\newtheorem{Proposition}[Theorem]{Proposition}

\begin{document}

\title{\textbf{Relative consistency of a finite nonclassical theory incorporating ZF and category theory with ZF}}

\author{Marcoen J.T.F. Cabbolet$^a$, Adrian R.D. Mathias$^b$\\
        $^a$\small{\textit{CLEA, Vrije Universiteit Brussel}}\footnote{email: Marcoen.Cabbolet@vub.be}\\
        $^b$\small{\textit{Department of Mathematics, Freiburg University}}\footnote{email: ardm@posteo.uk}
        }
\date{}

\maketitle

\begin{abstract}%\small
{ZF is undoubtedly the most widely accepted foundational theory for mathematics. However, ZF has two properties that already by the founding fathers of axiomatic set theory were considered pathological: ZF has infinitely many axioms, and ZF has a countable model.} Recently, in \emph{Axioms} \textbf{10}(2): 119 (2021), a nonclassical first-order theory $\mathfrak{T}$ of sets and functions has been introduced as the collection of axioms we have to accept if we want a foundational theory for (all of) mathematics that is not weaker than ZF but that lacks these two pathological properties. Here we prove that $\mathfrak{T}$ is relatively consistent with ZF. We conclude that this is a concrete step towards showing that $\mathfrak{T}$ is an advancement in the foundations of mathematics.
\end{abstract}

\section{Introduction}
The most widely accepted foundational theory for mathematics is without a doubt Zermelo-Fraenkel set theory. Its constructive axioms are those of Zermelo's set theory \cite{Zermelo}, plus the axiom schema of replacement suggested by Fraenkel in \cite{Fraenkel}. In accordance with common convention, we will use ZFC to denote the full theory, ZF for the full theory minus the axiom of choice (AC), and ZF(C) in statements that are to hold for both ZF and ZFC. It is then a fact that ZF(C) has infinitely many axioms. And as a corollary of the (downward) L\"{o}wenheim-Skolem theoren, developed by L\"{o}wenheim and Skolem \cite{Lowenheim,Skolem}, it is also a fact that if ZF(C) has a model, then it has a countable model. These two facts can be viewed as \emph{unwanted} or \emph{pathological features} of ZF(C).

Recently a finitely axiomatized nonstandard first-order theory of functions and sets, denoted by $\mathfrak{T}$, has been published as the collection of axioms one has to accept to get rid of these unwanted features \cite{Cabbolet}. The axioms of $\mathfrak{T}$ consist of some twenty standard first-order axioms, plus one nonstandard first-order axiom. That nonstandard axiom is so powerful that it admits the deduction of the infinite axiom schemata of separation and replacement (translated in the language of $\mathfrak{T}$) from the finitely many axioms of $\mathfrak{T}$. In addition, the (downward) L\"{o}wenheim-Skolem theorem does not hold for $\mathfrak{T}$: if $\mathfrak{T}$ has a model, then the model is uncountable.

{The purpose of this paper is to show that $\mathfrak{T}$ is relatively consistent with ZF. For that matter, this paper is organised as follows. For the sake of self-containment of this paper, the remainder of this introduction discusses the motivation for $\mathfrak{T}$ (Sect. 1.1) and summarily outlines $\mathfrak{T}$ (Sect. 1.2); the material is taken largely from \cite{Cabbolet}. It has to be kept in mind, however, that the summary outline of $\mathfrak{T}$ in Sect. 1.2 does not replace the formal introduction of $\mathfrak{T}$ in \cite{Cabbolet}. Sect. \ref{sect:method} outlines the method for the relative consistency proof. Sect. \ref{sect:result} presents the proof of relative consistency of $\mathfrak{T}$ with ZF. The final section discusses comparisons with other foundational approaches (Sects. \ref{sect:infinitary}--\ref{sect:TypeTheory}), and states the conclusions (Sect. \ref{sect:conclusion}).}

\subsection{{Motivation for the development of $\mathfrak{T}$}}
{With regards to the infinite axiomatization of ZF(C), of course one can shrug one's shoulders as a response. We note, however, that historically the development of ZF(C) has been followed by the development of Von Neumann-Bernays-G\"{o}del set theory (NBG) \cite{VonNeumann1928,VonNeumann1929,Bernays,Godel}, which can be finitely axiomatized, and by Montague's proof that ZF(C) cannot be finitely axiomatized \cite{Montague}. So, at the risk of committing the `\emph{post hoc ergo propter hoc}' fallacy, we conclude that these landmark developments must have been spurred by a dissatisfaction with the infinite axiomatization of ZF(C)---why else do all that work? We share that dissatisfaction: a foundation for mathematics should preferably be finitely axiomatized.}

{Then about the downward Loewenheim-Skolem theorem. It is not constructive, and as such the Loewenheim-Skolem theorem has no implications for everyday's mathematical practice. {But it implies that ZF(C) has a countable model: while this does not yield a logical inconsistency like Russel's paradox, it still is a pathological feature of ZF(C) because it yields the \emph{paradoxical state of affairs}---as Skolem called it---that the very axioms of ZF(C), which are intended to serve as a foundation for the existence of uncountable sets, can be satisfied by a merely countable model. Of course we can deny that this feature of ZF(C) is pathological by pointing out that the judgment `pathological' is made from a meta-theoretic perspective: \emph{internally} the power set of the natural numbers is uncountable, so we merely end up with some sort of ``relativity'' regarding the predicate `uncountable'. But this denial is then based on a rather isolationist view on mathematics that we do not agree with: our view is that mathematics provides a language for the sciences---\emph{mathematica ancilla scientiae}. In particular, when we do physics we work in a model of set theory. Now let us consider the case that we do physics in a countable model of ZF(C), and let us then give names to the elements $A\subseteq \mathbb{N}$ in $\mathcal{P}(\mathbb{N})$:
\begin{equation}\label{eq:name}
  A:= \sum\limits_{n=0}^{\infty} \chi(n) 2^{-n}
\end{equation}
where $\chi(n) = 0$ if $n\notin A$ and $\chi(n) = 1$ if $n\in A$. That way, the powerset $\mathcal{P}(\mathbb{N})$ becomes a subset of the real interval $[0,1]$. Regardless of whether we call this set `countable' from a meta-theoretic perspective or `uncountable' from an isolationist perspective, the point is that we definitely end up with a discrete geometry if we apply this set to model (segments of) physical dimensions. But there is a major problem with discrete geometries in physics. Let us consider linear motion of two particles on a plain, and let us first consider the view that the plane is continuous, that is, that the plane is $\mathbb{R}^2$ where $\mathbb{R}$ is the usual (uncountable) continuum. On the continuous plane, by setting off the particles at the right initial positions and with the right initial velocities, we can create the situation that the particles follow trajectories that intersect at a point $X$ on the plane, where they will collide. See Fig. \ref{fig:collisions}-(a) for an illustration. But now we consider the view that the plane is discrete: the particles then `leap' from point to point. But having been set off at the same initial positions and with the same initial velocities as before, the particles then do not collide when the previously mentioned point $X$ is not an element of the discrete space---see Fig. \ref{fig:collisions}-(b) for an illustration. These two views would give different predictions regarding the amount of annihilation events if we would apply them e.g. to crossing beams of electrons and positrons (who annihilate upon collision), with the discrete view giving the wrong prediction. So purely pragmatically, a countable model of ZF(C) is \emph{absolutely inadequate} for application to physics: the fact that ZF(C) has such a model is thus definitely a pathological feature from the perspective that \emph{mathematica ancilla scientiae}.}

\begin{figure}[h!]
\centering
\hfill
\subfigure[Motion in continuous space]{\includegraphics[width=0.425\textwidth]{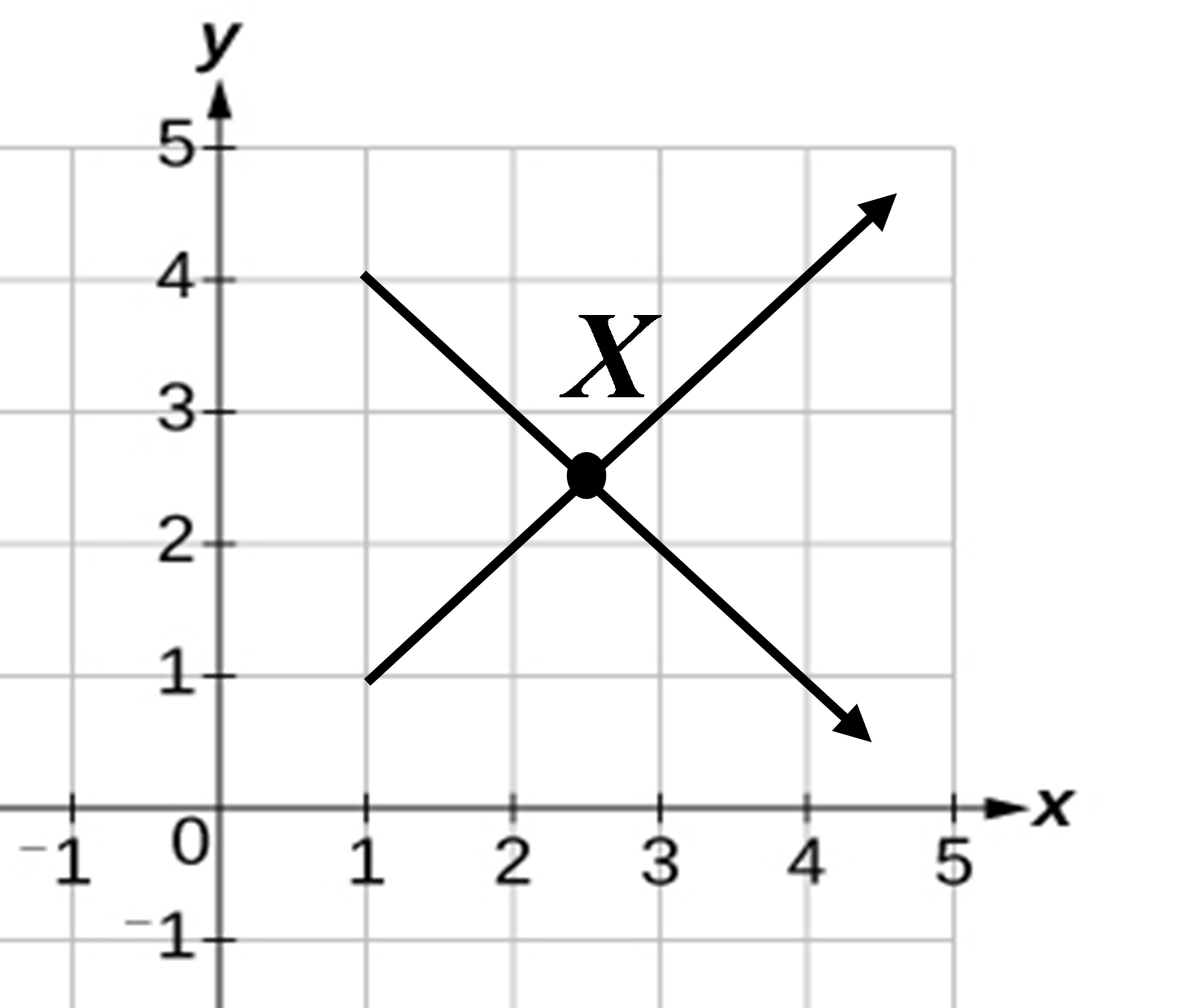}}
\hfill
\subfigure[Motion in discrete space]{\includegraphics[width=0.425\textwidth]{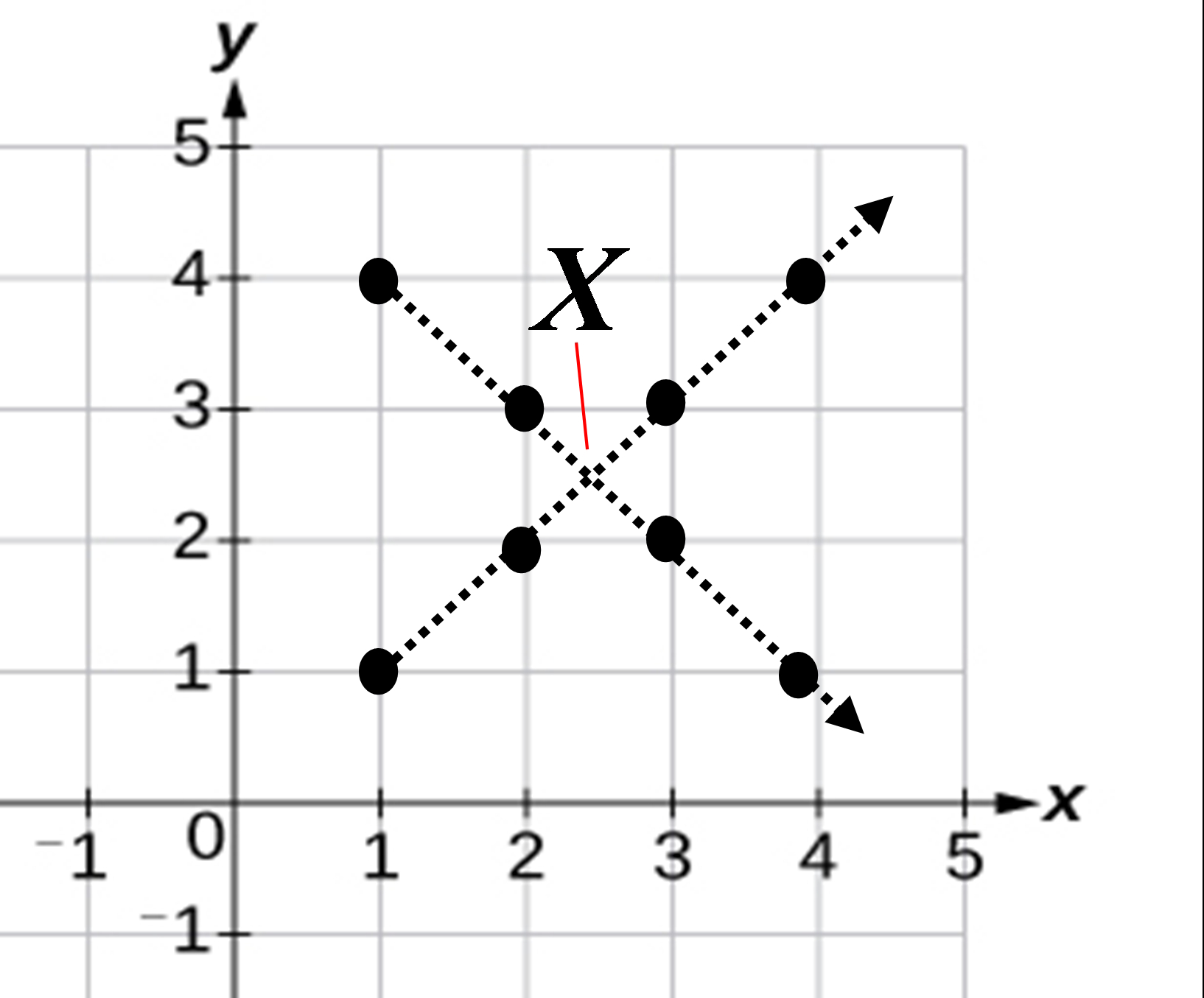}}
\hfill
%\subfigure[Diagram 3: four vertices]{\includegraphics[width=0.6\textwidth]{FeynmanDiagram3.jpg}}
%\hfill
\caption{{In diagram (a), two particles move linearly on continuous trajectories in $\mathbb{R}^2$ as indicated by the arrows: the particles collide at the point $X$ where the trajectories intersect. In diagram (b), the particles move with the same velocity in a discrete space, where they `leap' from point to point: the particles do not collide because the point $X$ where the dotted lines intersect is not an element of the discrete space.}}
\label{fig:collisions}
\end{figure}

{According to Bays, ``Skolem himself claimed that the paradox shows that set theory provides an inadequate foundation for mathematics'' \cite{Bays2}. However, Skolem's paper doesn't use a single, precise phrase to express this claim, nor was Skolem categorically against the idea of set theory as a foundation for mathematics in general: rather, his paper criticizes the idea of a set theory \emph{whose axioms are formulated in (standard) first-order language} as a foundation for mathematics by discussing the paradoxical nature of a countable model thereof. In his own words, as translated by Ebbinghaus in \cite{Ebbinghaus2},
\begin{quote}
  ``{\it I believed that it was so clear that axiomatization in terms of sets was not a satisfactory ultimate foundation of mathematics that mathematicians would, for the most part, not be very much concerned with it. But in recent times I have seen to my surprise that so many mathematicians think that these axioms of set theory provide the ideal foundation for mathematics; therefore it seemed to me that the time had come for a critique.}''
\end{quote}}
\noindent In addition, Zermelo himself of all people considered this countable model a pathological feature of ZF(C) \cite{Ebbinghaus}; he went on to develop second-order set theory, in which the theorem doesn't hold \cite{Zermelo2}. Now NBG shares this feature with ZF(C), and in that regard also Von Neumann has been quoted stating \cite{Ebbinghaus2}:
\begin{quote}
  ``{\it At present we can do no more than note that we have one more reason here to entertain reservations about set theory and that for the time being no way of rehabilitating this theory is known}''
\end{quote}
So with regard to this countable model of ZF(C), of course one can pretend to not care, or one can adopt a defeatist attitude that nothing can be done about it. But we go with Zermelo, Skolem, Von Neumann, and countless others: a foundational theory for mathematics should preferably lack this pathological property, even at the cost of departing from the framework of standard first-order logic. The Loewenheim-Skolem theorem has spurred considerable debate, which is beyond the scope of this paper; for an overview, see \cite{Bays}.}

{While getting rid of the two pathological properties of ZF has been the motivation for the development of $\mathfrak{T}$, the end result---i.e., the theory $\mathfrak{T}$---can be called \emph{stronger than} ZF. First of all, $\mathfrak{T}$ implies ZF, so all theorems that can be derived within ZF can also be derived within $\mathfrak{T}$. That establishes already that $\mathfrak{T}$ is \emph{not weaker than} ZF. But $\mathfrak{T}$ gets to this result with a lot less axioms---finitely many for $\mathfrak{T}$ versus infinitely many for ZF. Now in a weightlifting competition, if one and the same weight is lifted by two lifters, then the lifter with the lower body weight is considered \emph{stronger}. By analogy, $\mathfrak{T}$ is then stronger than ZF. But even if one doesn't like that argument, we still have the situation that we can prove a proposition in the framework of $\mathfrak{T}$ which we cannot prove in the framework of ZF, namely that the theory has no countable model: since this is a desirable result, we surely can call $\mathfrak{T}$ \emph{stronger than} ZF.}

\subsection{Summary outline of $\mathfrak{T}$}

%For the sake of self-containment of this paper, the nonstandard theory $\mathfrak{T}$ is summarily outlined in the remainder of this introduction. It has to be kept in mind, however, that this summary outline does not replace the formal introduction of $\mathfrak{T}$ in \cite{Cabbolet}.\\
%\ \\
To begin with, the universe of discourse is a category of sets and functions: its objects are \emph{sets}, its arrows are \emph{functions}. Importantly, while a function can be an element of a set, it is not a set: it is a concatenation of two sets. To explain, let $\{\alpha_1, \ldots\}, \{\alpha_2, \ldots\}, \ldots$ be any sets. Then a \textbf{concatenation} of $n$ sets is a thing
$$\{\alpha_1, \ldots\} \{\alpha_2, \ldots\} \cdots\{\alpha_n, \ldots\}$$
which for $n > 1$ is different from any set $Z$:
\begin{equation}\label{eq:concatenation}
  \{\alpha_1, \ldots\} \{\alpha_2, \ldots\} \cdots\{\alpha_n, \ldots\} \neq Z
\end{equation}
So, very simple, a concatenation of two sets is a thing $\{\alpha, \ldots\} \{\beta, \ldots\}$ which is not a set. (By the same token, an array of two bricks is not one brick.) A function is then a concatenation of two sets---the first is the graph of the function, the second is its domain---with the second one written as a subscript as in $F_X$. An \textbf{ur-function} in this framework is a function with a singleton domain and a singleton range. E.g. for Zermelo ordinals 0 and 1, the ur-function from $\{0\}$ to $\{1\}$ that maps 0 to 1 is the thing $\{ (0, 1)\}_{\{0\}}$ given by
\begin{gather}\label{eq:ur-function1}
\{ (0, 1)\}_{\{0\}}: \{0\}\twoheadrightarrow\{1\}\\
\{ (0, 1)\}_{\{0\}}: 0\mapsto 1 \label{eq:ur-function2}
\end{gather}
where the arrow $\twoheadrightarrow$ indicates a surjection. The ur-function $\{ (0, 1)\}_{\{0\}}$ can be graphically represented by a diagram; see Fig. \ref{fig:ur-function} below.
\vfill
\begin{SCfigure}[1.0][b!]
%\begin{center}
\includegraphics[width=0.5\textwidth]{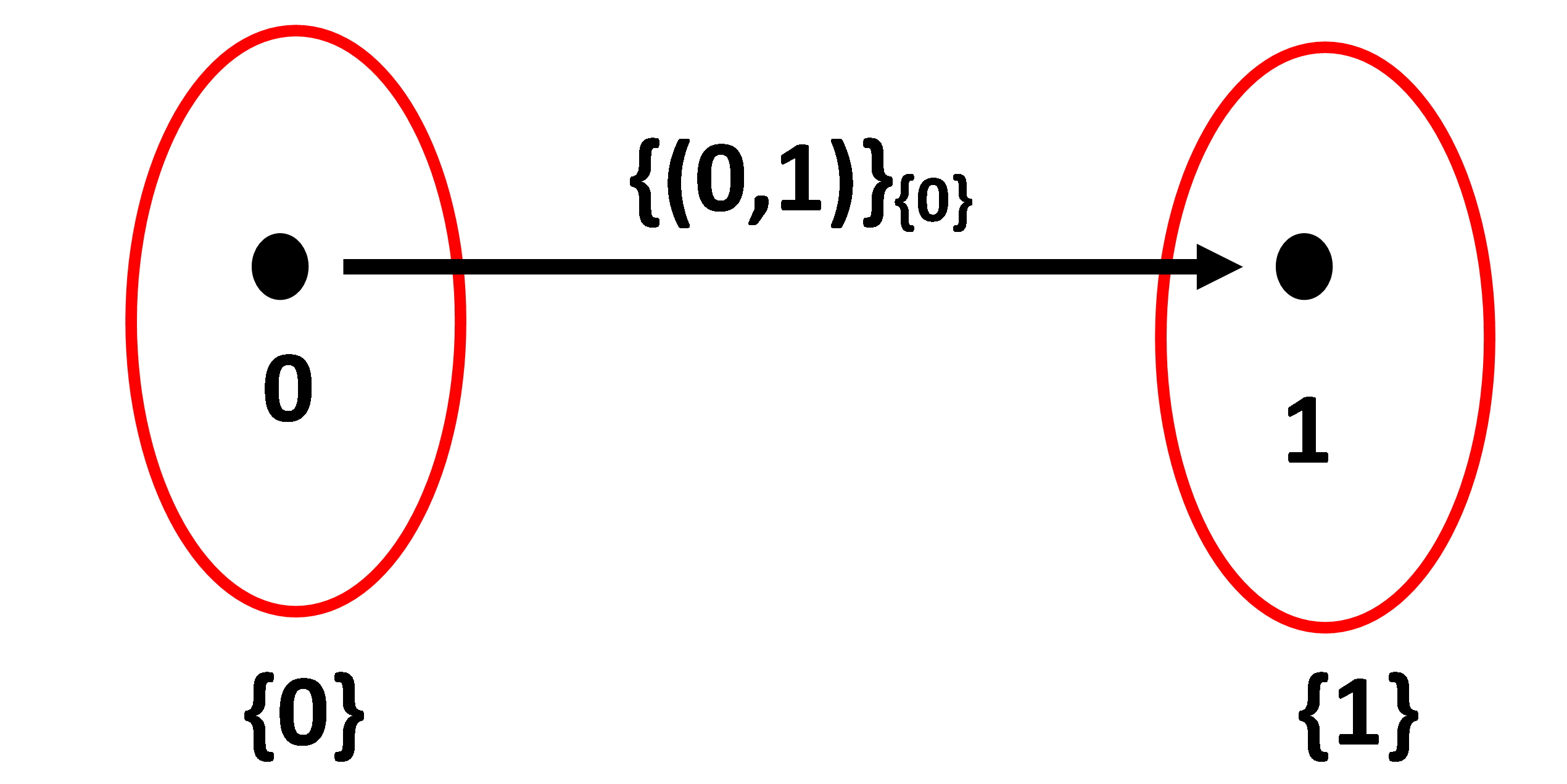}
%\end{center}
\caption{Diagram representing the ur-function $\{ (0, 1)\}_{\{0\}}$. On the left, a Venn diagram representing the singleton $\{0\}$. On the right, a Venn diagram representing the singleton $\{1\}$. The arrow represents the mapping $\{ (0, 1)\}_{\{0\}}: 0 \mapsto 1$.}
\label{fig:ur-function}
\end{SCfigure}
\newpage
The next thing is that all but one of the axioms of $\mathfrak{T}$ are expressed in a standard first-order language. So, we proceed by listing the terms of the language and its atomic expressions, starting with the terms:
\begin{itemize}
  \item constants $\emptyset$ and $\omega$ stand for the empty set and the (infinite) set of all finite Zermelo ordinals;
  \item (labeled) Roman variables $x, x_1, X, y, \ldots$ vary over sets;
  \item for a constant set $\mathbf{\hat{X}}$, variables $f_{\mathbf{\hat{X}}}, G_{\mathbf{\hat{X}}}, \ldots$ vary over functions on the set $\mathbf{\hat{X}}$;
  \item the constant $1_\emptyset$ stands for the inactive function;
  \item for a variable $X$ ranging over sets, composites variables $F_X, f_X, h_X, \ldots$ range over functions on a set $X$;
  \item (labeled) Greek variables $\alpha, \alpha_1, \beta, \ldots$ range over all things (set and functions).
\end{itemize}
The atomic expressions are the following:
\begin{itemize}
  \item $t_1 = t_2$, meaning $t_1$ is identical to $t_2$;
  \item $t_1 \in t_2$, meaning $t_1$ is an element of $t_2$;
  \item $t_1: t_2 \twoheadrightarrow t_3$, meaning $t_1$ is a surjection from $t_2$ onto $t_3$;
  \item $t_1: t_2 \mapsto t_3$, meaning $t_1$ maps $t_2$ to $t_3$.
\end{itemize}
The standard axioms of $\mathfrak{T}$ can then be divided into set-theoretical axioms and function-theoretical axioms. Without further explanation, the set-theoretical axioms are the following:
\begin{gather}
\forall X \forall Y (X = Y \Leftrightarrow \forall \alpha (\alpha \in X \Leftrightarrow \alpha \in Y))\\
\forall X \forall f_X \forall Y (f_X \neq Y)\\
\forall X \forall\alpha\forall\beta(X:\alpha\not\twoheadrightarrow\beta \wedge X:\alpha\not\mapsto\beta)\\
\exists X(X = \emptyset \wedge \forall\alpha(\alpha \not \in X))\\
\forall \alpha \forall \beta \exists X \forall \gamma (\gamma\in X \Leftrightarrow \gamma = \alpha \vee \gamma = \beta)\\
\forall X \exists Y \forall \alpha(\alpha \in Y \Leftrightarrow \exists Z(Z \in X \wedge \alpha \in Z))\\
\forall X \exists Y \forall \alpha (\alpha \in Y \Leftrightarrow \exists Z(Z \subset X \wedge \alpha = Z))\\
\forall X \neq \emptyset \exists \alpha (\alpha \in X \wedge \forall \beta (\beta \in \alpha \Rightarrow \beta \not \in X))\\
\exists X(X = \omega \wedge \emptyset \in \bm{\omega} \wedge \forall \alpha (\alpha\in \bm{\omega} \Rightarrow \alpha^+ \in \bm{\omega}) \wedge \forall \beta \in \bm{\omega} (\not\exists\gamma \in \bm{\omega}(\beta= \gamma^+)\Leftrightarrow\beta=\emptyset))
\end{gather}
This last axiom, the infinite set axiom, makes use of the function symbol $(.)^+$. This is defined as follows:
\begin{equation}
  \forall\alpha\forall\beta(\beta = \alpha^+ \Leftrightarrow \exists X(\beta = X \wedge \forall\gamma(\gamma\in X \Leftrightarrow\gamma = \alpha)))
\end{equation}
%\newpage
Proceeding, the standard function-theoretic axioms of $\mathfrak{T}$ are the following:
\begin{gather}
\forall X \forall f_X \forall \alpha (\alpha \not \in f_X)\\
\forall X \forall f_X (X\neq \emptyset \Rightarrow \exists Y\exists Z(f_X: Y \twoheadrightarrow Z \wedge \forall \alpha\in Y \exists! \beta (f_X: \alpha \mapsto \beta)))\\
\forall X \forall f_X \forall \alpha (\alpha\neq X \Rightarrow \forall \xi (f_X: \alpha \not \twoheadrightarrow \xi))\\
\forall X \forall f_X\forall \alpha \not\in X \forall \beta(f_X: \alpha \not\mapsto \beta)\\
\forall X \forall f_X (X\neq \emptyset \Rightarrow \forall \beta(f_X:X \twoheadrightarrow \beta \Rightarrow \exists Z(\beta = Z \wedge \forall \gamma(\gamma \in Z \Leftrightarrow \exists \eta\in X(f_X: \eta \mapsto \gamma)))))\\
\forall X \neq \emptyset \forall f_X\forall Y(f_X: X\twoheadrightarrow Y\Rightarrow  \forall \beta\exists Z\forall \alpha(\alpha \in Z \Leftrightarrow \alpha \in X \wedge f_X: \alpha \mapsto \beta))\\
\exists f_\emptyset (f_\emptyset = 1_\emptyset \wedge f_\emptyset: \emptyset \twoheadrightarrow \emptyset \wedge \forall \alpha \forall \beta (f_\emptyset: \alpha \not \mapsto \beta))\\
\forall \alpha \forall \beta \exists f_{\alpha^+} (f_{\alpha^+}: \alpha^+ \twoheadrightarrow \beta^+ \wedge f_{\alpha^+}: \alpha \mapsto \beta)\\
\forall X \forall f_X\forall Y(f_X: X\twoheadrightarrow Y\Rightarrow \forall \alpha(f_X: f_X \not\mapsto \alpha \wedge f_X:\alpha \not\mapsto f_X))
\end{gather}
With that, we have listed all standard axioms of $\mathfrak{T}$. But the power of the theory $\mathfrak{T}$ comes from its nonclassical Sum Function Axiom. We will demonstrate the axiom for a three-element set $S=\{\alpha,\beta,\gamma\}$ where $\alpha$, $\beta$, and $\gamma$ are any three different things (set or function), but for a rigorous introduction we refer the reader to the literature \cite{Cabbolet}. So, let any three ur-functions $f_{\{\alpha\}}$, $f_{\{\beta\}}$, and $f_{\{\gamma\}}$ be on the singletons of $\alpha$, $\beta$, and $\gamma$, respectively, be given; let the images of the ur-functions be denoted by $f_{\{\alpha\}}(\alpha)$, $f_{\{\beta\}}(\beta)$ and $f_{\{\gamma\}}(\gamma)$. Then the Sum Function Axiom guarantees that there is a sum function $F_S$ and a set $Y$ such that $F_S$ is a surjection from $S$ onto $Y$ and such that for each element of $S$, its image under sum function is the same as its image under the ur-function on the corresponding singleton. We can write this down in a first-order formula with composite variables $g_Z$ ranging over functions on the set $Z$:
\begin{equation}\label{eq:SumFunctionStandard}
%\begin{array}{l}
\forall f_{\{\alpha\}} \forall f_{\{\beta\}} \forall f_{\{\gamma\}}  \exists F_S\exists Y : %\\
F_S: S \twoheadrightarrow Y \wedge F_S:\alpha \mapsto f_{\{\alpha\}}(\alpha) \wedge F_S:\beta \mapsto f_{\{\beta\}}(\beta) \wedge F_S:\gamma \mapsto f_{\{\gamma\}}(\gamma)
%\end{array}
\end{equation}
This reads as: for any ur-function $f_{\{\alpha\}}$ on the singleton of $\alpha$ and for any ur-function $f_{\{\beta\}}$ on the singleton of $\beta$ and for any ur-function $f_{\{\gamma\}}$ on the singleton of $\gamma$ there is a function $F_{\{\alpha,\beta,\gamma\}}= F_S$ and a set $Y$ such that $F_S$ is a surjection from $S$ onto $Y$ and such that $F_S$ maps $\alpha$ to $f_{\{\alpha\}}(\alpha)$ and $\beta$ to $f_{\{\beta\}}(\beta)$ and $\gamma$ to $f_{\{\gamma\}}(\gamma)$. Using Venn diagrams, we can represent this graphically, as shown in Fig. \ref{fig:SumFunction} below.

\begin{figure}[b!]
\begin{center}
\includegraphics[width=0.95\textwidth]{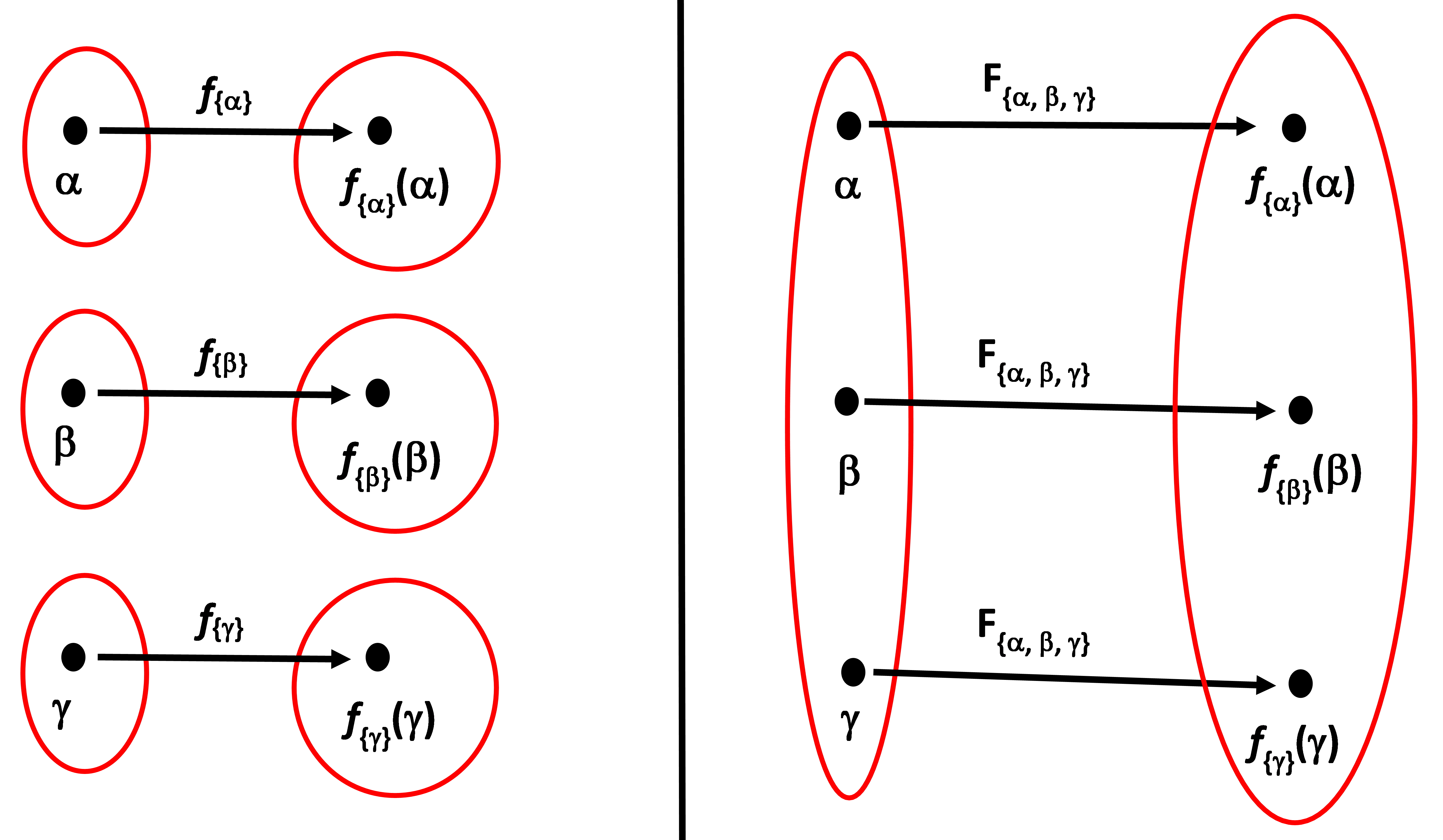}
\end{center}
\caption{Graphical representation of the sum function axiom for the three-element set $S=\{\alpha,\beta,\gamma\}$. The three diagrams on the left represent the three ur-functions $f_{\{\alpha\}}$, $f_{\{\beta\}}$, and $f_{\{\gamma\}}$ as in Fig. \ref{fig:ur-function}. The diagram on the right represents the postulated sum function $F_{\{\alpha,\beta,\gamma\}}$.}
\label{fig:SumFunction}
\end{figure}

Now we abbreviate formula \eqref{eq:SumFunctionStandard} by using a \textbf{multiple quantifier} $(\forall f_{\{\xi\}})_{\xi \in \{\alpha,\beta,\gamma\}}$ and a \textbf{conjunctive operator} $\bigwedge\limits_{\xi \in \{\alpha,\beta,\gamma\}}$ to an equivalent well-formed formula in the language of $\mathfrak{T}$:
\begin{equation}\label{eq:SumFunctionNonStandard}
(\forall f_{\{\xi\}})_{\xi \in \{\alpha,\beta,\gamma\}} \exists F_{\{\alpha,\beta,\gamma\}}\exists Y \left[
F_{\{\alpha,\beta,\gamma\}}: \{\alpha,\beta,\gamma\} \twoheadrightarrow Y \wedge
\bigwedge\limits_{\xi \in \{\alpha,\beta,\gamma\}} F_{\{\alpha,\beta,\gamma\}}:\xi \mapsto f_{\{\xi\}}(\xi)\right]
\end{equation}
This reads as: for any family of ur-functions $f_{\{\xi\}}$ indexed in $\{\alpha,\beta,\gamma\}$, there is a set $Y$ and a function $F_{\{\alpha,\beta,\gamma\}}$ such that $F_{\{\alpha,\beta,\gamma\}}$ is a surjection from $\{\alpha,\beta,\gamma\}$ onto the set $Y$ and, conjunctively over all $\xi \in \{\alpha,\beta,\gamma\}$, $F_{\{\alpha,\beta,\gamma\}}$ maps $\xi$ to its image under the ur-function $f_{\{\xi\}}$. It has to be taken that \emph{in this simple case} the multiple quantifier is equivalent to three standard quantifiers, and the subformula $\bigwedge\limits_{\xi \in \{\alpha,\beta,\gamma\}} F_{\{\alpha,\beta,\gamma\}}:\xi \mapsto f_{\{\xi\}}(\xi)$ is equivalent to a standard conjunction of three propositions---this may be obvious from comparing formulas \eqref{eq:SumFunctionStandard} and \eqref{eq:SumFunctionNonStandard}. The full Sum Function Axiom, abbreviated SUM-F, now says that formula \eqref{eq:SumFunctionNonStandard} holds for any (possibly infinite) non-empty set $X$---not just for the set $\{\alpha,\beta,\gamma\}$. So, SUM-F is the axiom
\begin{equation}\label{eq:Sum-F}
\forall X (\forall f_{\{\xi\}})_{\xi \in X}   \exists F_{X}\exists Y\left[
F_{X}: X \twoheadrightarrow Y \wedge
\bigwedge\limits_{\xi \in X} F_{X}:\xi \mapsto f_{\{\xi\}}(\xi)\right]
\end{equation}
This is a new nonstandard mathematical principle; in the general case, the multiple quantifier $(\forall f_{\{\xi\}})_{\xi \in X}$ can be equivalent to infinitely many standard quantifiers, and the subformula $\bigwedge\limits_{\xi \in X} F_{X}:\xi \mapsto f_{\{\xi\}}(\xi)$ can be equivalent to a standard conjunction of infinitely many propositions.
%\newpage

{The rules of inference are those of standard first-order logic, supplemented with a number of rather intuitive nonstandard rules that enable the deduction of standard sentences from SUM-F, plus a number of again rather intuitive nonstandard rules that enable the deduction of nonstandard sentences from standard sentences. Without further explanation, these nonstandard rules of inference are listed below; for a full treatment, see \cite{Cabbolet}.
\begin{enumerate}[(R-1)]
  \item Nonstandard Universal Elimination:
        {\setlength{\mathindent}{0cm}\begin{equation}
        \forall X\Psi \vdash [{\rm\hat{\mathbf{X}}}\backslash X]\Psi  \hfill {\rm for\ any\ constant\ {\rm\hat{\mathbf{X}}}}
        \end{equation}}
        where $\Psi$ is a formula with an occurrence of a nonstandard language elements.
  \item Multiple Universal Elimination:
        {\setlength{\mathindent}{0cm}\begin{equation}
        {(\forall f_{\alpha^+})}_{\alpha \in {\rm\hat{\mathbf{X}}}}\Phi \vdash [{\rm\hat{\mathbf{f}}}_{\alpha^+}\backslash f_{\alpha^+}]\Phi
        \end{equation}}
        where ${\rm\hat{\mathbf{f}}}_{\alpha^+}$ is a variable ranging over a specific family of ur-functions $({\rm\hat{\mathbf{f}}}_{\alpha^+})_{\alpha\in {\rm\hat{\mathbf{X}}}}$.
  \item Nonstandard Rule-C:
        {\setlength{\mathindent}{0cm}\begin{equation}
        \exists t \Phi \vdash [{\rm\hat{\mathbf{t}}}\backslash t]\Phi
        \end{equation}}
        where ${\rm\hat{\mathbf{t}}}$ is a constant in the range of $t$ that does not occur in $\Phi$ but for which $[{\rm\hat{\mathbf{t}}}\backslash t]\Phi$ holds.
  \item Conjunctive Operator Elimination:
        {\setlength{\mathindent}{0cm}\begin{equation}
        {\bigwedge}_{\alpha\in {\rm\hat{\mathbf{X}}}} \Psi(\alpha) \vdash [{\rm\bm{\hat{\alpha}}}\backslash \alpha] \Psi(\alpha)
        \end{equation}}
        where $\Psi(\alpha)$ is a formula of the type $t:t^\prime \mapsto t^{\prime\prime}$ that is open in $\alpha$, and ${\rm\bm{\hat{\alpha}}}$ any constant designating an element of ${\rm\hat{\mathbf{X}}}$.
  \item Conjunctive Operator Introduction:
        {\setlength{\mathindent}{0cm}\begin{equation}
        {\{ \Psi(\alpha) \}}_{\alpha\in \mathrm{\hat{\mathbf{X}}}}\vdash{\bigwedge}_{\alpha\in {\rm\hat{\mathbf{X}}}} \Psi(\alpha)
        \end{equation}}
        where $\Psi(\alpha)$ is an atomic formula of the type $t:t^\prime\mapsto t^{\prime\prime}$ that is open in $\alpha$.
  \item Multiple Universal Introduction:
        {\setlength{\mathindent}{0cm}\begin{equation}
        \Phi\left({{\bigwedge}_{\alpha\in {\rm\hat{\mathbf{X}}}}\Psi({\rm\hat{\mathbf{f}}}_{\alpha^+}) }\right)
        \vdash
        {(\forall f_{\alpha^+})}_{\alpha \in {\rm\hat{\mathbf{X}}}}[f_{\alpha^+}\backslash{\rm\hat{\mathbf{f}}}_{\alpha^+}]\Phi
        \end{equation}}
        where $\Phi\left({{\bigwedge}_{\alpha\in {\rm\hat{\mathbf{X}}}}\Psi({\rm\hat{\mathbf{f}}}_{\alpha^+}) }\right)$ denotes a formula $\Phi$ with a subformula ${\bigwedge}_{\alpha\in {\rm\hat{\mathbf{X}}}}\Psi({\rm\hat{\mathbf{f}}}_{\alpha^+})$ (implying that $\Psi$ is an atomic formula of the type $t:t^\prime\mapsto t^{\prime\prime}$), and where the variable ${\rm\hat{\mathbf{f}}}_{\alpha^+}$ ranges over an \textbf{arbitrary} family of ur-functions indexed in ${\rm\hat{\mathbf{X}}}$.
  \item Multiple Existential Introduction:
        {\setlength{\mathindent}{0cm}\begin{equation}
        \Phi\left({{\bigwedge}_{\alpha\in {\rm\hat{\mathbf{X}}}}\Psi({\rm\hat{\mathbf{f}}}_{\alpha^+}) }\right)
        \vdash
        {(\exists f_{\alpha^+})}_{\alpha \in {\rm\hat{\mathbf{X}}}}[f_{\alpha^+}\backslash{\rm\hat{\mathbf{f}}}_{\alpha^+}]\Phi
        \end{equation}}
        where $\Phi\left({{\bigwedge}_{\alpha\in {\rm\hat{\mathbf{X}}}}\Psi({\rm\hat{\mathbf{f}}}_{\alpha^+}) }\right)$ denotes a formula $\Phi$ with a subformula ${\bigwedge}_{\alpha\in {\rm\hat{\mathbf{X}}}}\Psi({\rm\hat{\mathbf{f}}}_{\alpha^+})$ (implying that $\Psi$ is an atomic formula of the type $t:t^\prime\mapsto t^{\prime\prime}$), and where the variable ${\rm\hat{\mathbf{f}}}_{\alpha^+}$ ranges over a \textbf{specific} family of ur-functions indexed in ${\rm\hat{\mathbf{X}}}$.
\end{enumerate}
This works as follows. If we have constructed a set designated by the constant ${\rm\hat{\mathbf{X}}}$, e.g. the set of natural numbers $\mathbb{N}$, then using rule (R-1) we can deduce from SUM-F that
\begin{equation}
{(\forall f_{\alpha^+})}_{\alpha \in {\rm\hat{\mathbf{X}}}} \exists F_{\rm\hat{\mathbf{X}}} \exists Y\left(F_{\rm\hat{\mathbf{X}}}:{\rm\hat{\mathbf{X}}} \twoheadrightarrow Y \wedge {\bigwedge}_{\alpha\in {\rm\hat{\mathbf{X}}}} F_{\rm\hat{\mathbf{X}}}:\alpha\mapsto f_{\alpha^+}(\alpha)\right)
\end{equation}
If we then have constructed a family of ur-functions indexed in ${\rm\hat{\mathbf{X}}}$, e.g. the family of ur-functions $({\rm\hat{\mathbf{f}}}_{\alpha^+})_{\alpha\in \mathbb{N}}$ given by ${\rm\hat{\mathbf{f}}}_{\alpha^+}: \alpha \mapsto \alpha+1$ for any $\alpha\in \mathbb{N}$, then using rule (R-2) we can deduce that
\begin{equation}
\exists F_{\rm\hat{\mathbf{X}}} \exists Y\left(F_{\rm\hat{\mathbf{X}}}:{\rm\hat{\mathbf{X}}} \twoheadrightarrow Y \wedge {\bigwedge}_{\alpha\in {\rm\hat{\mathbf{X}}}} F_{\rm\hat{\mathbf{X}}}:\alpha\mapsto {\rm\hat{\mathbf{f}}}_{\alpha^+}(\alpha)\right)
\end{equation}
From there we deduce, using rule (R-3), that
\begin{equation}\label{eq:RuleC}
\exists Y({\rm\hat{\mathbf{F}}}_{\rm\hat{\mathbf{X}}}:{\rm\hat{\mathbf{X}}}\twoheadrightarrow Y) \wedge
{\bigwedge}_{\alpha\in {\rm\hat{\mathbf{X}}}} {\rm\hat{\mathbf{F}}}_{\rm\hat{\mathbf{X}}}:\alpha\mapsto ({\rm\hat{\mathbf{f}}}_{\alpha^+}(\alpha)
\end{equation}
which uses the constant ${\rm\hat{\mathbf{F}}}_{\rm\hat{\mathbf{X}}}$ designating the sum function on ${\rm\hat{\mathbf{X}}}$. So in the example with $\mathbb{N}$, this sum function is the function on $\mathbb{N}$ that maps any natural number $\alpha$ to $\alpha+1$, as follows from using (R-4).}
%\newpage

{This may seem cumbersome to use in everyday's mathematical practice, but luckily we don't have to. It is, namely, a theorem of $\mathfrak{T}$---actually, it is a theorem schema---that for any nonempty set $X$, if there is a functional relation $\Phi(\alpha, \beta)$ that relates every $\alpha$ in $X$ to precisely one $\beta$, then there is a function $F_X$ with some codomain $Y$ that maps every $\eta \in X$ to precisely that $\xi \in Y$ for which $\Phi(\eta, \xi)$. In a formula, using the iota-operator:
\begin{gather}
\forall X\neq\emptyset (\forall \alpha \in X \exists! \beta \Phi(\alpha, \beta) \Rightarrow \exists F_X\exists Y (F_X:X\twoheadrightarrow Y \wedge \forall \eta \in X(F_X: \eta \mapsto \imath\xi\Phi(\eta, \xi))))
\end{gather}
This is the main theorem of $\mathfrak{T}$, which makes $\mathfrak{T}$ easy to use in everyday mathematical practice because for any set $X$ we can construct a function $f_X$ by giving a defining function prescription $f_X:X\twoheadrightarrow Y\ , \ f_X: \alpha\overset{\rm def}{\longmapsto}\iota\beta\Phi(\alpha,\beta)$ where $\Phi$ is some functional relation: $\mathfrak{T}$ then guarantees that $F_X$ exists, as well as its graph, its image set, and the inverse image sets for every element of its codomain---ergo, giving a defining function prescription is a tool for \mbox{constructing sets.}}

{But not just that. The main theorem also allows us to derive the theorem schema of separation, which is a translation of the corresponding axiom schema of ZF in the language of $\mathfrak{T}$:
\begin{gather}
  \forall X \exists Y \forall \alpha : \alpha \in Y \Leftrightarrow \alpha \in X \wedge \Phi(\alpha)
\end{gather}
with subformulas $\Phi$ and $\Psi$ interpreted as usual. So, working in the framework of $\mathfrak{T}$, one can use separation just as one would in the framework of ZF.}

\section{Method}\label{sect:method}
We want to prove that $\mathfrak{T}$ is relative consistent with ZF. That is, we want to prove that
\begin{equation}\label{eq:TrcwZF}
  \left[ \nvdash_{ZF} \bot\right] \Rightarrow \left[ \nvdash_{\mathfrak{T}} \bot\right]
\end{equation}
We have to note, however, that the symbol `$\vdash$' on the left hand side in Eq. \eqref{eq:TrcwZF} represents standard inference, while the symbol `$\vdash$' on the right hand side in Eq. \eqref{eq:TrcwZF} represents inference in the nonclassical logical framework of $\mathfrak{T}$. Thus speaking, if the symbol `$\Vvdash$' represents inference in the nonstandard logical framework of $\mathfrak{T}$, then actually the following proposition has to be proved:
\begin{equation}\label{eq:TrcwZF2}
  \left[ \not\vdash_{ZF} \bot\right] \Rightarrow \left[ \not\Vvdash_{\mathfrak{T}} \bot\right]
\end{equation}
We will prove this proposition in three steps. In the first step, we define the standard first-order theory $\mathfrak{T}^\infty_0$ as the theory $\mathfrak{T}$ minus the nonstandard Sum Function Axiom but extended with a collection of standard first-order theorems of $\mathfrak{T}$ (Def. \ref{def:Tinfty0}), and we prove that $\mathfrak{T}^\infty_0$ is relative consistent with ZF as in
      \begin{equation}\label{eq:T0rcwZF}
      \left[ \nvdash_{ZF} \bot\right] \Rightarrow \left[ \nvdash_{\mathfrak{T}^\infty_0} \bot\right]
      \end{equation}
{In this first step, we use Set Matrix Theory (SMT), introduced in \cite{SMT} as a generalisation of ZF in which finite matrices with set-valued entries are treated as \emph{objects sui generis}; an outline of SMT is given in the Appendix. The motivation for its use is purely pragmatical: in \cite{CabboletRC} it has been proved that SMT is relative consistent with ZF, and therefore the use of SMT provides a fairly easy way to prove Eq. \eqref{eq:T0rcwZF}. Namely, if SMT is relatively consistent with ZF, then so is a subtheory SMT$_{1\times2}$ of SMT in which only $1\times2$ set matrices $[G\ \ D]$ occur whose entries $G$ and $D$ can be viewed as the graph and the domain of a function, and so is a conservative extension SMT$_{1\times2}^{+4}$ of SMT$_{1\times2}$ obtained by adding definitions for expressions $t_1: t_2 \twoheadrightarrow t_3$ and $t_1: t_2 \mapsto t_3$ and for constants $\omega$ and $[\emptyset\ \ \emptyset]$---we then have
\begin{equation}
\left[ \nvdash_{ZF} \bot\right] \Rightarrow \left[ \nvdash_{\rm SMT_{1\times2}^{+4}} \bot\right]
\end{equation}
And when $G$ and $D$ are graph and domain of a function, there is not really any difference between an object $[G\ \ D]$ in the language of SMT$_{1\times2}^{+4}$ and an object $G_D$ in the language of $\mathfrak{T}^\infty_0$, so by proving that there is an interpretation of ${\mathfrak{T}^\infty_0}$ in SMT$_{1\times2}^{+4}$ we prove
\begin{equation}
\left[ \nvdash_{\rm SMT_{1\times2}^{+4}} \bot\right] \Rightarrow \left[ \nvdash_{\mathfrak{T}^\infty_0} \bot\right]
\end{equation}
Eq. \eqref{eq:T0rcwZF} then follows by transitivity of relative consistency.} In the second step, we show that $\mathfrak{T}^\infty_0$ in the nonstandard logical framework of $\mathfrak{T}$ is relatively consistent with $\mathfrak{T}^\infty_0$ in the framework of standard first-order logic, as in
      \begin{equation}\label{eq:DeSwart}
      \left[ \not\vdash_{\mathfrak{T}^\infty_0} \bot\right] \Rightarrow \left[ \not\Vvdash_{\mathfrak{T}^\infty_0} \bot\right]
      \end{equation}
Finally, in the third step we prove that $\mathfrak{T}$ is relatively consistent with $\mathfrak{T}^\infty_0$ as in
  \begin{equation}\label{eq:TrcwT0}
  \left[ \not\Vvdash_{\mathfrak{T}^\infty_0} \bot\right] \Rightarrow \left[ \not\Vvdash_{\mathfrak{T}} \bot\right]
  \end{equation}
By transitivity of relative consistency, we may then infer that $\mathfrak{T}$ is relatively consistent with ZF as in Eq. \eqref{eq:TrcwZF2}. The next three sections treat the above three steps one by one.

\section{Relative consistency of $\mathfrak{T}$ with ZF}\label{sect:result}
\subsection{Relative consistency of $\mathfrak{T}^\infty_0$ with ZF}
Starting from Set Matrix Theory (SMT), we proceed as follows:
\begin{enumerate}[(i)]
  \item we define STM$_{1\times2}$ as a restriction of SMT in which only set matrices of the type $[f\ \ X]$ occur where $f$ is the graph of a function on a set $X$ (Def. \ref{def:SMT1x2});
  \item we show that that SMT$_{1\times2}$ is relatively consistent with SMT (Prop. \ref{prop:SMT1x2RCwSMT});
  \item we define SMT$^{+4}_{1\times2}$ as an extension of STM$_{1\times2}$ with definitions for expressions $t_1: t_2 \twoheadrightarrow t_3$ and $t_1: t_2 \mapsto t_3$ and for constants $\omega$ and $[\emptyset\ \ \emptyset]$ (Def. \ref{def:SMT1x2+4});
  \item we show that SMT$^{+4}_{1\times2}$ is relatively consistent with STM$_{1\times2}$ (Prop. \ref{prop:SMT1x2+5RCwSMT});
  \item we show that $\mathfrak{T}^\infty_0$ is relatively consistent with SMT$^{+4}_{1\times2}$ (Prop. \ref{prop:Tinfty0=SMT-}).
\end{enumerate}
We may then infer that $\mathfrak{T}^\infty_0$ is relatively consistent with SMT. By transitivity of relative consistency, we may then infer that $\mathfrak{T}^\infty_0$ is relatively consistent with ZF.\\
\ \\
To begin with, recall from \cite{Cabbolet} that the theorem schemata of separation and replacement of $\mathfrak{T}$ are the schemata
\begin{gather}\label{eq:SeparationScheme}
  \forall X \exists Y \forall \alpha : \alpha \in Y \Leftrightarrow \alpha \in X \wedge \Phi(\alpha)\\
\label{eq:SubstitutionScheme}
  \forall X : \forall \alpha \in X \exists! \beta \Psi(\alpha, \beta) \Rightarrow \exists Z \forall \gamma( \gamma \in Z \Leftrightarrow \exists \xi : \xi \in X \wedge \Psi(\xi, \gamma))
\end{gather}
with subformulas $\Phi$ and $\Psi$ interpreted as usual. In addition, the following is a theorem of $\mathfrak{T}$:
\begin{equation}\label{eq:ReverseGraphTheorem}
\forall f\forall X\forall Y : f\in Y^X \Rightarrow \exists h_X\forall\alpha\forall\beta (h_X:\alpha\mapsto\beta \Leftrightarrow (\alpha,\beta)\in f)
\end{equation}
where $Y^X$ is the set of all graphs of functions from $X$ to $Y$; the proof is omitted.
\begin{Definition}\label{def:Tinfty0}
Let $\Gamma$ be the collection of formulas in the language of $\mathfrak{T}$ made up of all theorems in the schemata \eqref{eq:SeparationScheme} and \eqref{eq:SubstitutionScheme} plus theorem \eqref{eq:ReverseGraphTheorem}. Let $\mathfrak{T}[\Gamma]$ be the supertheory of $\mathfrak{T}$ obtained by adding the formulas of $\Gamma$ as axioms to $\mathfrak{T}$. Then $\mathfrak{T}^\infty_0$ is the subtheory of $\mathfrak{T}[\Gamma]$ obtained by removing the nonstandard sum function axiom from the axioms of $\mathfrak{T}[\Gamma]$, and by removing the nonstandard deduction rules. \hfill$\blacksquare$
\end{Definition}

\noindent The theory $\mathfrak{T}^\infty_0$ has, thus, infinitely many axioms (as indicated by the superscript $\infty$) but zero nonstandard axioms (as indicated by the subscript $0$). The axioms of $\mathfrak{T}^\infty_0$ are, thus, the standard set-theoretical and function-theoretical axioms of $\mathfrak{T}$ listed in the introduction, plus all formulas in $\Gamma$ as listed above. %In order to show that $\mathfrak{T}^\infty_0$ is relatively consistent with SMT$^{+4}_{1\times2}$, we first have to define SMT$^{+4}_{1\times2}$.

\begin{Definition}\label{def:SMT1x2}
Let the theory SMT$_{1\times2}$ be the subtheory of SMT obtained by
\begin{enumerate}[(i)]
  \item removing all function symbols $f_{m\times n}$ from the language, except $f_{1\times2}$---note that this eliminates the reduction axiom and the omission axiom scheme, given by Eqs. \eqref{eq:reduction} and \eqref{eq:omission} in the Appendix \ref{app:SMT}---so that Greek variables $\alpha,\beta, \ldots$ now range over sets and $1\times2$ matrices;
  \item replacing the set matrix axiom schema, given by Eq. \eqref{eq:SetMatrixAxiom} in Appendix \ref{app:SMT}, by
        {\setlength{\mathindent}{0cm}\begin{equation}\label{eq:1x2axiom}
        \left(\forall f \forall X \forall Y : f\in Y^X \Rightarrow \exists\alpha (\alpha = [f \ \ X])\right) \wedge \left(\forall \alpha : \exists\beta\exists\gamma(\alpha =[\beta\ \ \gamma]) \Rightarrow \exists X\exists f\exists Y (\alpha = [f\ \ X] \wedge f\in Y^X)\right)
        \end{equation}}
  \item removing the set of set matrices axiom scheme, given by Eq. \eqref{eq:SetOfSetMatricesAxiom} in Appendix \ref{app:SMT}
\end{enumerate}
$\blacksquare$
\end{Definition}
\newpage
\noindent In the framework of SMT we can construct set matrices $\left [ \begin{array}{ccc} \alpha_{11} & \ldots & \alpha_{1n} \\ \vdots &   & \vdots \\ \alpha_{m1} & \ldots & \alpha_{mn} \end{array} \right ]$ of any size $m\times n$ and with any entries $\alpha_{ij}$. In the framework of SMT$_{1\times2}$, however, we only want to allow $1\times2$ set matrices $[f\ \ X]$ with set-valued entries $f$ and $X$ satisfying the condition that $f$ is the graph of a function on $X$. This restriction forces us to remove the set of set matrices axiom scheme, because that would introduce the existence of $1\times2$ set matrices $[\alpha\ \ \beta]$ violating this condition. Altogether, we then end up with Def. \ref{def:SMT1x2}.

As a side note, in the framework of SMT, a cartesian product of sets $A$ and $B$ is just the set of $1\times2$ matrices $[\alpha\ \ \beta]$ with $\alpha\in A$ and $\beta\in B$, which can be constructed from the set of set matrices $M_{1\times2}(A\cup B)$ and a separation axiom. But in the framework of SMT$_{1\times2}$ we can no longer construct the cartesian product $A\times B$ that way, because we removed the set of set matrices axiom scheme. So in order to be able construct the cartesian product $A\times B$ in the framework of SMT$_{1\times2}$, we are forced to define a two-tuple $\langle \alpha, \beta\rangle$ as a set---we can use any one of the usual definitions, e.g. $\langle\alpha, \beta\rangle = \{\alpha, \{\alpha, \beta\}\}$.
%\newpage
\begin{Proposition}\label{prop:SMT1x2RCwSMT}
SMT$_{1\times2}$ is relatively consistent with SMT.\hfill$\blacksquare$
\end{Proposition}

\paragraph{Proof:} The proof is omitted.\hfill$\blacksquare$

\begin{Definition}\label{def:SMT1x2+4}
The theory SMT$^{+4}_{1\times2}$ is the theory SMT$_{1\times2}$ extended with
\begin{enumerate}[(i)]
  \item a definition of a ternary relation $t_1: t_2 \twoheadrightarrow t_3$:
        {\setlength{\mathindent}{0cm}\begin{equation}\label{eq:surjection}
        \forall\alpha\forall\beta\forall\gamma \left( \alpha:\beta\twoheadrightarrow\gamma\Leftrightarrow\exists f\exists X\exists Y : \alpha = [f\ \ X] \wedge \beta = X \wedge \gamma = Y \wedge \forall \mu \in Y\exists \nu \in X(\langle \nu, \mu\rangle \in f) \right)
        \end{equation}}
  \item a definition of a ternary relation $t_1: t_2 \mapsto t_3$:
        {\setlength{\mathindent}{0cm}\begin{equation}\label{eq:mapping}
        \forall\alpha\forall\beta\forall\gamma \left( \alpha:\beta\mapsto\gamma\Leftrightarrow\exists X\exists f : \alpha = [f\ \ X] \wedge \langle \beta, \gamma\rangle \in f \right)
        \end{equation}}
  \item an axiom for the constant $\omega$ representing the set of all Zermelo ordinals:
        {\setlength{\mathindent}{0cm}\begin{equation}\label{eq:omega}
        \forall x: x = \omega \Leftrightarrow \emptyset \in x \wedge \forall \alpha (\alpha\in x \Rightarrow \alpha^+ \in x) \wedge \forall \beta \in x (\not\exists\gamma \in x(\beta= \gamma^+)\Leftrightarrow\beta=\emptyset))
        \end{equation}}
  \item an axiom for the constant $[\emptyset\ \ \emptyset]$:
        {\setlength{\mathindent}{0cm}\begin{equation}\label{eq:InactiveFunction}
        \forall \alpha: \alpha = [\emptyset\ \ \emptyset] \Leftrightarrow \alpha:\emptyset\twoheadrightarrow\emptyset \wedge \forall\beta\forall\gamma (\alpha:\beta\not\mapsto \gamma)
        \end{equation}}
\end{enumerate}
$\blacksquare$
\end{Definition}
%\newpage
\noindent So in words, in the framework of SMT$^{+4}_{1\times2}$ we have that the object $\alpha$ is a surjection from $\beta$ to $\gamma$ if and only if there are sets $f$, $X$, and $Y$ such that
\begin{enumerate}[(i)]
  \item $\alpha$ is the $1\times2$ set matrix $[f\ \ X]$---note that by Eq. \eqref{eq:1x2axiom}, this implies that $f$ is the graph of a function on $X$;
  \item $\beta$ is the set $X$;
  \item $\gamma$ is the set $Y$ and for every element $\mu$ in $Y$ there is an element $\nu$ in $X$ such that $\langle \nu, \mu\rangle$ is in $f$.
\end{enumerate}
Furthermore, the object $\alpha$ maps $\beta$ to $\gamma$ if and only if there is a set $X$ and a graph $f$ of a function on $X$ such that $\alpha$ is the $1\times2$ set matrix $[f\ \ X]$ and the two-tuple $\langle \beta, \gamma\rangle$ is in $f$.

\begin{Proposition}\label{prop:SMT1x2+5RCwSMT}
SMT$^{+4}_{1\times2}$ is relatively consistent with SMT.\hfill$\blacksquare$
\end{Proposition}

\paragraph{Proof:} It suffices to show that SMT$^{+4}_{1\times2}$ is a conservative extension of SMT$_{1\times2}$: in that case SMT$^{+4}_{1\times2}$ is relatively consistent with SMT$_{1\times2}$, and given Prop. \ref{prop:SMT1x2RCwSMT} we may then infer the proposition.

Firstly, Eqs. \eqref{eq:surjection} and \eqref{eq:mapping} are \emph{defining axioms} of the ternary predicates that generate the atomic expressions $t_1: t_2 \twoheadrightarrow t_3$ and $t_1: t_2 \mapsto t_3$. The extension of SMT$_{1\times2}$ by Eqs. \eqref{eq:surjection} and \eqref{eq:mapping} is therefore an extension by definitions: this is a conservative extension \cite{Shoenfield}.

Secondly, Eqs. \eqref{eq:omega} and \eqref{eq:InactiveFunction} are \emph{defining axioms} of the constants $\omega$ and $[\emptyset\ \ \emptyset]$. The further extension of SMT$_{1\times2}$ by \eqref{eq:omega} and \eqref{eq:InactiveFunction} is therefore again an extension by definitions. Ergo, SMT$^{+4}_{1\times2}$ is a conservative extension of SMT$_{1\times2}$.\hfill$\blacksquare$
\newpage
\begin{Proposition}\label{prop:Tinfty0=SMT-}
$\mathfrak{T}^\infty_0$ is relatively consistent with SMT$^{+4}_{1\times2}$.\hfill$\blacksquare$
\end{Proposition}

\paragraph{Proof:} It suffices to define an interpretation of $\mathfrak{T}^\infty_0$ in SMT$^{+4}_{1\times2}$, such that the interpretation of every axiom of $\mathfrak{T}^\infty_0$ is a theorem of SMT$^{+4}_{1\times2}$ \cite{Shoenfield}. We define the interpretation function $\tau$ from the language $L_{\mathfrak{T}^\infty_0}$ of $\mathfrak{T}^\infty_0$ to the language $L_{{\rm SMT}^{+4}_{1\times2}}$ of SMT$^{+4}_{1\times2}$ as follows:
\begin{itemize}
  \item $\tau: \emptyset \mapsto \emptyset$ , $\tau: \omega \mapsto \omega$ , $\tau: 1_\emptyset\mapsto [\emptyset \ \ \emptyset]$;
  \item a simple variable $x$ ranging over sets is mapped by $\tau$ to a simple variable $x$ ranging over sets;
  \item a simple variable $\alpha$ ranging over all things (sets and functions) is mapped by $\tau$ to a simple variable $\alpha$ ranging over all things (sets and set matrices);
  \item for any constant $\bm{\hat{X}}$ of $L_{\mathfrak{T}^\infty_0}$ referring to an individual set, a composite term $f_{\bm{\hat{X}}}$ ranging over functions on the set $\bm{\hat{X}}$ is mapped by $\tau$ to a $1\times2$ matrix $[f\ \ \bm{\hat{X}}]$ which effectively ranges over functions on the set $\bm{\hat{X}}$;
  \item a composite variable $f_X$ ranging over functions on a set $X$ is mapped by $\tau$ to a $1\times2$ matrix $[f\ \ X]$ which effectively ranges over functions on a set $X$;
  \item atomic formulas $t_1 \in t_2$, $t_1 = t_2$, $t_1: t_2 \twoheadrightarrow t_3$, and $t_1: t_2 \mapsto t_3$ are mapped by $\tau$ to respectively formulas $\tau(t_1) \in \tau(t_2)$, $\tau(t_1) = \tau(t_2)$, $\tau(t_1): \tau(t_2)\twoheadrightarrow \tau(t_3)$, and $\tau(t_1): \tau(t_2)\mapsto \tau(t_3)$;
  \item formulas $\neg\Phi$ and $\Phi\ldots\Psi$ with logical connective `$\ldots$' are mapped by $\tau$ to formulas $\neg\tau(\Phi)$ and $\tau(\Phi)\ldots\tau(\Psi)$;
  \item for any simple variable $\alpha$ ranging over all things, a formula $\forall\alpha:\Phi$ is mapped by $\tau$ to a formula $\forall\alpha:\tau(\Phi)$;
  \item for any simple variable $X$ ranging over sets, a formula $\forall X:\Phi$ is mapped by $\tau$ to a formula $\forall X:\tau(\Phi)$;
  \item for any constant $\bm{\hat{X}}$ of $L_{\mathfrak{T}^\infty_0}$ referring to an individual set, a formula $\forall f_{\bm{\hat{X}}}:\Phi$ or $\exists f_{\bm{\hat{X}}}:\Phi$ is mapped by $\tau$ to a formula $\forall f:\exists\alpha\left(\alpha = [f\ \ \bm{\hat{X}}]\right)\Rightarrow\tau(\Phi)$ or $\exists f: \exists\alpha\left(\alpha = [f\ \ \bm{\hat{X}}]\right)\wedge\tau(\Phi)$;
  \item for any simple variable $X$ ranging over sets, a formula $\forall f_X:\Phi$ or $\exists f_X :\Phi$ in the scope of the quantifier $\forall X$ is mapped by $\tau$ to a formula $\forall f:\exists\alpha\left(\alpha = [f\ \ X]\right)\Rightarrow \tau(\Phi)$ or $\exists f: \exists\alpha\left(\alpha = [f\ \ X]\right)\wedge \tau(\Phi)$
\end{itemize}
All axioms of $\mathfrak{T}^\infty_0$ are then interpreted as theorems of SMT$^{+4}_{1\times2}$. We will prove this for the formulas of the infinite collection $\Gamma$, which are all axioms of $\mathfrak{T}^\infty_0$; the rest is left as an exercise for the reader. So, for starters, observe that the axiom schemata \eqref{eq:SeparationScheme} and \eqref{eq:SubstitutionScheme} of $\Gamma$ are interpreted \emph{identically} as the separation axiom scheme \eqref{eq:SEP} and the replacement axiom scheme \eqref{eq:GeneralizedREP} listed in Appendix \ref{app:SMT}: these are axiom schemata of SMT$^{+4}_{1\times2}$, with the Greek variables ranging over sets and $1\times2$ matrices $[f\ \ X]$ with $f$ being the graph of a function on $X$. Next, axiom \eqref{eq:ReverseGraphTheorem} of $\Gamma$ is interpreted as
\begin{equation}\label{eq:translation}
\forall f\forall X\forall Y : f\in Y^X \Rightarrow \exists h (\exists\gamma(\gamma = [h\ \ X]) \wedge \forall\alpha\forall\beta ([h\ \ X]:\alpha\mapsto\beta \Leftrightarrow (\alpha,\beta)\in f)
\end{equation}
To prove that this is a theorem of SMT$^{+4}_{1\times2}$, we assume
\begin{equation}\label{eq:proof1}
f\in Y^X
\end{equation}
for arbitrary $f, X, Y$. From assumption \eqref{eq:proof1} and axiom \eqref{eq:1x2axiom} we infer
\begin{equation}\label{eq:proof2}
\exists\gamma: \gamma = [f\ \ X]
\end{equation}
From Eq. \eqref{eq:proof2} and axiom \eqref{eq:mapping} we then infer
\begin{equation}
\forall\alpha\forall\beta ([f\ \ X]:\alpha\mapsto\beta \Leftrightarrow (\alpha,\beta)\in f)
\end{equation}
From here we deduce
\begin{equation}
f\in Y^X \Rightarrow \exists\gamma: \gamma = [f\ \ X]) \wedge \forall\alpha\forall\beta ([f\ \ X]:\alpha\mapsto\beta \Leftrightarrow (\alpha,\beta)\in f)
\end{equation}
We then obtain Eq. \eqref{eq:translation} by introducing existential and universal quantifiers.\hfill $\blacksquare$

\begin{Corollary}
$\mathfrak{T}^\infty_0$ is relatively consistent with ZF.\hfill$\blacksquare$
\end{Corollary}

\subsection{Extending the framework with nonclassical rules of inference}

To prove Eq. \eqref{eq:DeSwart}, we add the nonstandard rules of inference to the standard logical framework of $\mathfrak{T}^\infty_0$, and show that this doesn't introduce any inconsistencies. For that matter, it suffices to prove the following proposition.

\begin{Proposition}\label{prop:LogicalExtension}
Let $\Sigma$ be a collection of standard first-order formulas in the language $L(\mathfrak{T})$ of the nonstandard theory $\mathfrak{T}$. Then, if no contradiction can be derived from $\Sigma$ using the standard rules of inference of first-order logic, no contradiction can be derived from $\Sigma$ using the standard and nonstandard rules of inference of the axiomatic system of $\mathfrak{T}$:
\begin{equation}
\left[\Sigma \not \vdash \bot \right] \Rightarrow \left[\Sigma \not \Vvdash \bot \right]
\end{equation}
$\blacksquare$
\end{Proposition}

\noindent \textbf{Proof}: Let $\Sigma_1$ be the collection of formulas made up of the formulas $\phi\in\Sigma$ and the formulas that can be derived from formulas in $\Sigma$ by \textbf{\underline{once}} applying the rules of inference of standard first-order logic. Let $\Sigma_2$ be the collection of formulas made up of the formulas $\phi\in\Sigma_1$ and the formulas that can be derived from formulas in $\Sigma_1$ by \textbf{\underline{once}} applying the rules of inference of standard first-order logic, and so forth. Let $\Sigma_\omega$ be the collection of formulas such that $\phi \in \Sigma_\omega$ if and only if there is an integer $n$ for which $\phi\in\Sigma_n$. We now assume
\begin{equation}
\Sigma \not \vdash \bot
\end{equation}
So, there is no formula $\phi$ such that $\phi\in\Sigma_\omega \wedge \neg\phi\in\Sigma_\omega$. Now there are just a handful of nonstandard rules of inference in the axiomatic system of $\mathfrak{T}$, so we start by treating the first one: introduction of a conjunctive operator $\bigwedge\limits_{\alpha\in \mathbf{\hat{X}}}$ for a constant $\mathbf{\hat{X}}\in L(\mathfrak{T})$, given by
\begin{equation}
\left\{ \mathbf{\hat{F}}_{\mathbf{\hat{X}}} : \alpha\mapsto \mathbf{\hat{f}}_{\{\alpha\}}(\alpha) \right\}_{\alpha\in \mathbf{\hat{X}}}  \Vvdash  \bigwedge\limits_{\alpha \in \mathbf{\hat{X}}} \mathbf{\hat{F}}_{\mathbf{\hat{X}}}:\alpha \mapsto \mathbf{\hat{f}}_{\{\alpha\}}(\alpha)
  \end{equation}
This is rule R-5 of Sect. 1.2. Now let $\Upsilon_1$ be the collection of formulas made up of the formulas $\phi\in\Sigma_\omega$ and the formulas that can be derived from formulas in $\Sigma_\omega$ by \textbf{\underline{once}} applying rule R-5. Since rule R-5 only applies to collections of formulas of the type $t_1: t_2\mapsto t_3$, this only adds formulas of the type $\bigwedge\limits_{\alpha \in \mathbf{\hat{X}}} \mathbf{\hat{F}}_{\mathbf{\hat{X}}}:\alpha \mapsto \mathbf{\hat{f}}_{\{\alpha\}}(\alpha)$ to $\Sigma_\omega$. But such a formula cannot be the negation of a formula $\phi\in\Sigma_\omega$, because it contains a new logical symbol (conjunctive operator) that does not occur in any of the formulas that are in $\Sigma_\omega$. Ergo, there is no formula $\phi$ such that $\phi\in\Upsilon_1 \wedge \neg\phi\in\Upsilon_1$.

Proceeding, let $\Upsilon_2$ be the collection of formulas made up of the formulas $\phi\in\Upsilon_1$ and the formulas that can be derived from formulas in $\Upsilon_1$ by
\textbf{\underline{once}} applying one of the standard rules or rule R-5, and so forth; let $\Upsilon_\omega$ be the collection of formulas such that $\phi \in \Upsilon_\omega$ if and only if there is an integer $n$ for which $\phi\in\Upsilon_n$. But rule R-5 only applies to collections of atomic formulas of the type $t_1: t_2\mapsto t_3$---it doesn't apply to any collection of propositions as in infinitary logic---and that means that we cannot create new nonstandard formulas of the type $\bigwedge\limits_{\alpha \in \mathbf{\hat{X}}} \mathbf{\hat{F}}_{\mathbf{\hat{X}}}:\alpha \mapsto \mathbf{\hat{f}}_{\{\alpha\}}(\alpha)$ when we go from $\Upsilon_n$ to $\Upsilon_{n+1}$. And that means that $\Upsilon_\omega$ consists of all the formulas that can be derived from the formulas $\phi\in\Upsilon_1$ by standard propositional logic. Ergo, given the consistency of standard propositional logic, there is no formula $\phi$ such that $\phi\in\Upsilon_\omega \wedge \neg\phi\in\Upsilon_\omega$.

We now add the second nonstandard rule of inference: elimination of a conjunctive operator $\bigwedge\limits_{\alpha\in \mathbf{\hat{X}}}$ for constants $\bm{\hat{\alpha}}, \mathbf{\hat{X}}\in L(\mathfrak{T})$, with $\bm{\hat{\alpha}}\in \mathbf{\hat{X}}$, given by
\begin{equation}
\bigwedge\limits_{\alpha \in \mathbf{\hat{X}}} \mathbf{\hat{F}}_{\mathbf{\hat{X}}}:\alpha \mapsto \mathbf{\hat{f}}_{\{\alpha\}}(\alpha) \Vvdash \mathbf{\hat{F}}_{\mathbf{\hat{X}}}:\bm{\hat{\alpha}} \mapsto \mathbf{\hat{f}}_{\{\bm{\hat{\alpha}}\}}(\bm{\hat{\alpha}})
  \end{equation}
This is rule R-4 of Sect. 1.2. Let $\Phi_1$ then be the collection of formulas made up of the formulas $\phi\in\Upsilon_\omega$ and the formulas that can be derived from formulas in $\Upsilon_\omega$ by \textbf{\underline{once}} applying rule $I_2$. But by applying rule R-4, we can only derive formulas $\phi$ of the type $t_1: t_2\mapsto t_3$ that were already in $\Sigma_\omega$. So, if we define $\Phi_{n+1}$ as the collection of formulas made up of the formulas $\phi\in\Phi_n$ and the formulas that can be derived from formulas in $\Phi_n$ by \textbf{\underline{once}} applying one of the standard rules or one of the nonstandard rules R-4 or R-5, and we define $\Phi_\omega$ as the collection of formulas such that $\phi \in \Phi_\omega$ if and only if there is an integer $n$ for which $\phi\in\Phi_n$, then there is no formula $\phi$ such that $\phi\in\Phi_\omega \wedge \neg\phi\in\Phi_\omega$.

The other nonstandard inference rules (introduction of a multiple quantifier, elimination of a multiple quantifier, introduction of a single quantifier for a nonstandard formula, elimination of a single quantifier for a nonstandard formula) can be treated the same way---this is quite laborious, but straightforward. This gives a collection of formulas $\Omega_\omega$ that can be derived from the collection $\Sigma$ by applying the standard and nonstandard rules of inference of the axiomatic system of $\mathfrak{T}$: there is then no formula $\phi$ such that $\phi\in\Omega_\omega \wedge \neg\phi\in\Omega_\omega$. That means that $\Sigma$ is also consistent in the axiomatic system of $\mathfrak{T}$: we then have $\Sigma\not\Vvdash\bot$. Prop. \ref{prop:LogicalExtension} then follows immediately. \hfill$\blacksquare$

\begin{Corollary}
\begin{equation}\label{cor:LogicalExt}
\left[\mathfrak{T}^\infty_{0} \not \vdash \bot \right] \Rightarrow \left[\mathfrak{T}^\infty_{0} \not \Vvdash \bot \right]
\end{equation}
$\blacksquare$
\end{Corollary}

\subsection{Relative consistency of $\mathfrak{T}$ with $\mathfrak{T}^\infty_0$}\label{sect:NonstandardPart}

\noindent Our proof method is as follows. We assume $\mathfrak{T}^\infty_{0} \not \Vvdash \bot $. That is, we assume that $\mathfrak{T}^\infty_{0}$ is consistent. If we then prove that the negation of SUM-F cannot be derived from $\mathfrak{T}^\infty_{0}$ (Prop. \ref{prop:NotNotSum-F}) we may derive that $\mathfrak{T}^\infty_{0} \cup \{{\rm SUM-F}\}$ is consistent. But $\mathfrak{T}^\infty_{0} \cup \{{\rm SUM-F}\}$ is just the theory $\mathfrak{T}[\Gamma]$ of Def. \ref{def:Tinfty0}, which is, of course, equivalent to $\mathfrak{T}$. Ergo, if we prove Prop. \ref{prop:NotNotSum-F}, then we may conclude that $\mathfrak{T}$ is relatively consistent with $\mathfrak{T}^\infty_{0}$ as in Eq. \eqref{eq:TrcwT0}.

\begin{Proposition}\label{prop:NotNotSum-F}
The negation of SUM-F cannot be derived from $\mathfrak{T}^\infty_{0}$:
\begin{equation}
\mathfrak{T}^\infty_{0} \not\Vvdash \neg {\rm SUM-F}
\end{equation}
$\blacksquare$
\end{Proposition}

\paragraph{Proof:} SUM-F is the nonstandard axiom of $\mathfrak{T}$ given by Eq. \eqref{eq:Sum-F}; its negation is the formula
\begin{equation}\label{eq:NotSum-F}
\exists X (\exists f_{\{\xi\}})_{\xi \in X} : \neg\exists F_{X} \exists Y \left(F_{X}: X \twoheadrightarrow Y
\wedge \bigwedge\limits_{\xi \in X} F_{X}:\xi \mapsto f_{\{\xi\}}(\xi)\right)
\end{equation}
So, we want to prove that it cannot be deduced from $\mathfrak{T}^\infty_{0}$ that there is a set $X$ and a family of ur-functions $f_{\{\xi\}}$ indexed in $X$ such that there is no function $F_X$ and no set $Y$ such that $F_X$ is a surjection from $X$ onto $Y$ and, conjunctively over $\xi\in X$, $F_X$ maps $x$ to its image under the corresponding ur-function $f_{\{\xi\}}$. We proceed by \emph{reductio ad absurdum}: we assume Eq. \eqref{eq:NotSum-F}. Applying Rule-C, i.e. the rule $\exists \alpha\Phi(\alpha)\vdash \Phi(\bm{\hat{\alpha}})$ where $\bm{\hat{\alpha}}$ is a constant for which $\Psi$ is true, that means that there is a certain individual set $\mathbf{\hat{X}}$ and a certain family of ur-functions $(\mathbf{\hat{f}}_{\{\xi\}})_{\xi \in \mathbf{\hat{X}}}$ such that
\begin{equation}\label{eq:ToDisprove}
\neg\exists F_{\mathbf{\hat{X}}} \exists Y \left(F_{\mathbf{\hat{X}}}: \mathbf{\hat{X}} \twoheadrightarrow Y \right)
\wedge \bigwedge\limits_{\xi \in \mathbf{\hat{X}}} F_{\mathbf{\hat{X}}}:\xi \mapsto \mathbf{\hat{f}}_{\{\xi\}}(\xi)
\end{equation}
However, we have that
\begin{gather}
\forall \xi\in \mathbf{\hat{X}}\exists!\alpha : \alpha = \mathbf{\hat{f}}_{\{\xi\}}(\xi)\\
\forall \xi\in \mathbf{\hat{X}}\exists!\beta : \beta = (\xi, \mathbf{\hat{f}}_{\{\xi\}}(\xi))
\end{gather}
By the axiom schema of substitution \eqref{eq:SubstitutionScheme} of $\mathfrak{T}^\infty_{0}$, we thus obtain that
\begin{gather}
\exists Y: Y = \{\mathbf{\hat{f}}_{\{\xi\}}(\xi) \ | \ \xi \in \mathbf{\hat{X}}\}  \\
\exists G: G = \{(\xi, \mathbf{\hat{f}}_{\{\xi\}}(\xi)) \ | \ \xi \in \mathbf{\hat{X}}\}
\end{gather}
meaning that there are sets $Y$ and $G$ such that $Y$ is the set of all images of elements of $\mathbf{\hat{X}}$ under the ur-functions $(\mathbf{\hat{f}}_{\{\xi\}})_{\xi \in \mathbf{\hat{X}}}$ and $G$ is the set that contains a two-tuple $(\xi, \mathbf{\hat{f}}_{\{\xi\}}(\xi))$ for every $\xi\in \mathbf{\hat{X}}$. But then $G$ is the graph of a function on $\mathbf{\hat{X}}$, so by axiom \eqref{eq:ReverseGraphTheorem} of $\mathfrak{T}^\infty_{0}$ we have
\begin{equation}
\exists F_{\mathbf{\hat{X}}} \forall \xi\in \mathbf{\hat{X}} \left( F_{\mathbf{\hat{X}}} : \xi\mapsto \mathbf{\hat{f}}_{\{\xi\}}(\xi) \right)
\end{equation}
Using the nonstandard deduction rule $ \left\{ \mathbf{\hat{F}}_{\mathbf{\hat{X}}} : \xi\mapsto \mathbf{\hat{f}}_{\{\xi\}}(\xi) \right\}_{\xi\in \mathbf{\hat{X}}}  \Vvdash  \bigwedge\limits_{\xi \in \mathbf{\hat{X}}} \mathbf{\hat{F}}_{\mathbf{\hat{X}}}:\xi \mapsto \mathbf{\hat{f}}_{\{\xi\}}(\xi)$ we derive from here that
\begin{equation}
\exists F_{\mathbf{\hat{X}}} \exists Y \left(F_{\mathbf{\hat{X}}}: \mathbf{\hat{X}} \twoheadrightarrow Y
\wedge \bigwedge\limits_{\xi \in \mathbf{\hat{X}}} F_{\mathbf{\hat{X}}}:\xi \mapsto \mathbf{\hat{f}}_{\{\xi\}}(\xi)\right)
\end{equation}
which contradicts \eqref{eq:ToDisprove}. Ergo, our assumption \eqref{eq:NotSum-F} cannot be true. That proves Prop. \ref{prop:NotNotSum-F}.\hfill$\blacksquare$

\section{Discussion and conclusion}
{Mathematicians are not the same: what one mathematician thinks is acceptable as a foundational theory for mathematics will not be the opinion of all other mathematicians. Many find category theory a good setting for their research, others prefer a set-theoretic setting. In \cite{Cabbolet} a setting has been proposed which includes both of these, thus promoting the unity of mathematics. The sections 4.1 - 4.6 below compare this to other approaches to a foundation for mathematics. Sect. 4.7 states the conclusion.}
\subsection{Comparison to traditional infinitary frameworks}\label{sect:infinitary}

{Traditional infinitary frameworks have never gained any traction beyond the foundational debate. The main disadvantage is that the syntax of the language allows conjunctions of infinitely many formulas, without any further restrictions. Consequently, the language contains infinitary conjunctions that cannot be written down. But the primary function of the language of mathematics is \emph{written communication}, so this is a real argument against traditional infinitary frameworks. In the literature, such infinitary conjunctions are often referred to by designators like `$\bigwedge\limits_{i<\lambda}\phi_i$', where $\lambda$ is some infinite ordinal. The thing is, however, that the designator `$\bigwedge\limits_{i<\lambda}\phi_i$' is not a well-formed formula of the language: the well-formed formula to which it refers cannot be written down. E.g. the conjunction of all axioms of the separation axiom schema of ZF would be a single well-formed formula, but it cannot be written down. So, it is true that ZF can be finitely axiomatized in an infinitary framework, but then its axioms cannot be written down.}

{One might object to the claim that written communication is the primary function of the language of mathematics that it downplays the role that abstraction and idealization play in modern mathematics. In reply, however, we would like to point out that abstraction and idealization are processes intended to frame ideas in the language of mathematics, but in the end the outcomes of these processes are to be communicated. After all, science is a social enterprise aimed at the development of knowledge by an open discussion of new ideas. Mathematics provides a language for this open discussion, which is universal across languages and cultures, and which is concise and unambiguous---unlike natural language, which can be vague. Furthermore, this open discussion of new ideas takes place by means of \emph{written} communication. Many mathematically formulated ideas are too abstract or complex to be clearly communicated verbally, especially without visual aids. Written symbols allow for the layering of meaning, like nested functions, integrals, or logical statements, which are hard to follow by ear alone. Writing allows enormous compression of information. For instance, the compact notation $\eta_{\alpha\beta} x^\alpha y ^\beta$ in tensor calculus, assuming the Einstein summation convention, replaces a long and clunky verbal enumeration. But not just that. Written mathematics allows for careful structuring of arguments: theorems, lemmas, proofs, diagrams, etc. This structure is essential for rigor and verification---others can check each line of a reasoning, which is difficult to do in spoken form alone. Moreover, science builds on previous work. Written records allow results to be referenced, corrected, expanded, or disproved: without written form, complex results would be lost or distorted over time. That may justify the claim that written communication is the primary function of the language of mathematics.}

Compared to traditional infinitary frameworks, the syntax of the language of $\mathfrak{T}$ is much more restricted. Suppose that the constant $\mathbf{\hat{X}}$ refers to an infinite set; then for a variable $\alpha$, the conjunctive operator `$\bigwedge\limits_{\alpha\in\mathbf{\hat{X}}}$' can only be put in front of a formula $\Phi(\alpha)$ of the type $t_1:t_2\mapsto t_3$ that is open in $\alpha$. For example, consider the infinite set $\bm{\omega} = \{0,1,2,\ldots\}$ and let the function $\mathbf{\hat{F}}_{\bm{\omega}}$ be defined by $\mathbf{\hat{F}}_{\bm{\omega}}:\alpha \mapsto\alpha+1$. Then the expression $\bigwedge\limits_{\alpha\in\bm{\omega}} \mathbf{\hat{F}}_{\bm{\omega}}:\alpha \mapsto\alpha+1$, which (semantically) is equivalent to an infinitary conjunction
\begin{equation}\label{eq:InfinitaryConjunction}
\mathbf{\hat{F}}_{\bm{\omega}}:0 \mapsto 1 \wedge \mathbf{\hat{F}}_{\bm{\omega}}:1 \mapsto2 \wedge \mathbf{\hat{F}}_{\bm{\omega}}:2\mapsto3 \wedge \ldots
\end{equation}
in a traditional infinitary framework, actually is a well-formed formula in the language of $\mathfrak{T}$. But contrary to the infinitary conjunction \eqref{eq:InfinitaryConjunction}, the expression `$\bigwedge\limits_{\alpha\in\bm{\omega}} \mathbf{\hat{F}}_{\bm{\omega}}:\alpha \mapsto\alpha+1$' has a finite number of symbols: it can actually be written down.

{So, it is true that the syntax of the language of $\mathfrak{T}$ allows the construction of well-formed formulas that are semantically equivalent to infinitary conjunctions of infinitely many standard formulas, but these well-formed formulas are all finite: that is the main distinction between $\mathfrak{T}$ and traditional infinitary frameworks. Consequently, the main advantage of $\mathfrak{T}$ over traditional infinitary frameworks is that the criticism against traditional infinitary frameworks, being that the formulas cannot be written down, does not apply to $\mathfrak{T}$.}

\subsection{{Comparison to mainstream set theories}}
{Von Neumann-Bernays-G\"{o}del set theory (NBG), developed after ZF, has been proved to be a conservative extension of ZFC \cite{Felgner}. But NBG can be finitely axiomatized \cite{Godel}, and as such it lacks the unwanted property of ZF of having an infinite axiomatization. However, NBG shares with ZF the second pathological property of having a countable model.}

{So, we can prefer $\mathfrak{T}$ over ZF because (i) all theorems of ZF are also theorems of $\mathfrak{T}$, and (ii) $\mathfrak{T}$ lacks the two pathological properties of ZF. And we can prefer $\mathfrak{T}$ over NBG (i) all theorems of NBG that involve sets are also theorems of $\mathfrak{T}$, and (ii) $\mathfrak{T}$ lacks the pathological property of NBG of having a countable model. By virtue of the fact that all theorems of ZF are also theorems of $\mathfrak{T}$, one can do mathematics in the framework of $\mathfrak{T}$ just as one is used to do in the framework of ZF, yet without having to face the ugly truth that there is a countable model in which the powerset of the set of natural numbers is countable.}%\\

\subsection{{Comparison to constructive set theories}}

{Constructivist approaches to axiomatic set theory have led to a multiple constructive set theories, such as
\begin{itemize}
  \item intuitionistic ZF \cite{Friedman};
  \item constructive set theory \cite{Myhill,Friedman2};
  \item constructive ZF \cite{Aczel}.
\end{itemize}
The common denominator is that these all use the language of ZF in the framework of intuitionistic logic. While the various versions of constructive set theory differ with respect to their axioms, the number of axioms is infinite in all cases. Consequently, constructive set theories share with ZF the property of having an infinite axiomatization. Therefore, the present theory $\mathfrak{T}$ is preferable over constructive set theories because it has finitely many axioms. It would be interesting to find out whether $\mathfrak{T}$ can be recontextualized in the framework of intuitionistic logic, and if so, what remains of its power. That is left as a topic for further research.}

\subsection{{Comparison to the category-theoretical approach to set theory}}
{In 1965, Lawvere published his elementary theory of the category of sets (ETCS) \cite{Lawvere}. The intended universe of discourse is, at first glance, similar to that of $\mathfrak{T}$: a category, whose objects are sets and whose arrows are functions. The language, however, is that of category theory: in the framework of ETCS, the expression `$x\in A$', the most recognizable expression of set theory, is not an atomic expression of the language. Instead, it is formalized in the language of category theory as
\begin{equation}\label{eq:element}
  x\in A \Leftrightarrow x: 1\rightarrow A
\end{equation}
where 1 is a singleton that is assumed to exist.

ETCS, however, lacks an axiom of regularity. {As already mentioned by Von Neumann \cite{VonNeumann1925}, the universe of ETCS can therefore contain the following pathological sets:
\begin{itemize}
  \item the narcissus set $N$, which is full of itself as in $N = \{N\}$;
  \item Siamese twin sets $A$ and $B$, for which $A\in B\wedge B\in A$;
  \item sets giving rise to a Mirimanoff chain $x_1\in x_2\in\ldots\in x_n\in x_1$.
\end{itemize}
For a discussion of these sets and a proof how the axiom of regularity prevents their existence, see \cite{Muller}. What makes these sets `pathological' is that they are not well-founded, giving rise to infinitely descending $\in$-chains: these sets cannot be constructed with the constructive axioms of set theory, but without an axiom of regularity they can exist a priori in the mathematical universe. As Von Neumann put it, it would be desirable to get rid of these superfluous sets \cite{VonNeumann1925}.} Furthermore, ETCS is weaker than ZF: only ETCS[REP], i.e. ETCS extended with an axiom schema of replacement, is as strong as ZFC in the sense that the same theorems hold \cite{Leistner}.}
%\newpage

The present theory $\mathfrak{T}$ is then preferable over ETCS because ETCS is weaker and allows pathological objects. And the present theory $\mathfrak{T}$ is preferable over ETCS[REP] because $\mathfrak{T}$ is finitely axiomatized whereas ETCS[REP] shares with ZFC the property of having infinitely many axioms. In addition, ETCS[REP] still lacks the axiom of regularity: at thus also allows pathological objects that cannot exist in the universe of $\mathfrak{T}$.

An additional argument is the language: the replacement of `$x\in A$' by `$x:1\rightarrow A$' is not an improvement. {Namely, for every set $X$ containing the empty set, we now get an expression $\emptyset:1\rightarrow X$. Of course we can assume this formally as an axiom: it's not wrong in the sense that it yields a logical contradiction. But it gives the idea that the empty set is a function, which is conceptually very hard to reconcile with the idea of the empty set as just a set which does not contain any elements. In other words, the replacement of $x\in A$ by $x:1\rightarrow A$ does not yield an improvement in conceptual clarity---quite to the contrary, actually. On that ground, we don't agree with this way of reducing the language of mathematics to the language of category theory. Instead, the language of $\mathfrak{T}$ contains the language of set theory and the language of category theory.}

\subsection{{Comparison to countabilism}}

{`Countabilism' is the view that all sets are countable \cite{Barton}. So, if we view set theories that imply countabilism, then the Loewenheim-Skolem theorem is not problematic or pathological. And assuming we can get to a finite axiomatization, this too seems a viable approach to get to a foundational theory for mathematics that lacks the two pathological problems of ZF.}

{However, since mathematics provides the language for the sciences, the problem with countabilism lies in its application to physics: if we accept a set theory that implies countabilism as the foundation of mathematics, we run into seemingly unsolvable problems in physics because countable sets have to be used to represent real dimensions.
{Builes and Wilson have attempted to defend the position that physics is unproblematic in the framework of a mathematical foundation that implies countabilism \cite{Builes}:
\begin{quote}
``{\it the Countabilist has a number of options. First, they could reformulate orthodox theories using a regions-based conception of the continuum as developed by Hellman and Shapiro (2018). Second, the Countabilist (along with everyone else) could point out that our best physical theories cannot be strictly speaking true anyway.}''
\end{quote}
But this defense of countabilism in physics is inadequate. Looking at this regions-based conception of the continuum developed by Hellman and Shapiro in \cite{Hellman}, the basic idea is that a spatial dimension can be modeled by a set $D$ given by
\begin{equation}\label{eq:D}
D = \{(a, b)\ |\ a\ {\rm and}\ b\ {\rm are\ finite\ sums\ of\ integer\ powers\ of\ 2}\} \subset \{(x, y)\ |\ x,y\in \mathbb{Q}\}
\end{equation}
(The intervals considered by Hellman and Shapiro should be viewed as ur-elements rather than as uncountable sets of points, but nevertheless Eq. \eqref{eq:D} reflects the basic idea.) Viewed on its own, this is a self-consistent philosophy of space. The time dimension then also has to be modeled by a set like $D$. But it is, then, absolutely not clear what it means that a point-particle (like a photon, which doesn't have a spatial extension) has a ``position'' in space at a ``point'' in time. Worse, it is absolutely not clear what ``motion'' is in this framework, and how it avoids the picture in Fig. \ref{fig:collisions}-(b) which is certain to be wrong. So, it is not a certainty that orthodox theories of physics can \emph{at all} be reformulated using this region-based idea of the continuum.}

{A second problem with countabilism in physics is that it is not at all clear what a differential $dx$ is in a spatial dimension represented by a set like $D$ in Eq. \eqref{eq:D}. To defend countabilism in physics, Hellman and Shapiro refer to Arntzenius, who concludes from the fact that quantum mechanics (QM) uses Hilbert spaces with a countable basis that it is more natural to interpret wave functions as living on a pointless space \cite{Arntzenius}. The point (pun intended) is that in orthodox QM, these wave functions are actually complex functions on $\mathbb{R}^3$ (with $\mathbb{R}$ again being the uncountable set of the reals). But if we want to reformulate orthodox QM along the lines suggested by Hellman and Shapiro in the above quote, then the domain of the wave functions becomes a space $D^3$, but what is then their range? The complex numbers are out of the window with countabilism. Even more so, for a quantum system with one spatial dimension, its wave function $\psi$ evolves in time according to the Schroedinger equation
\begin{equation}\label{eq:SE}
i\hbar\frac{d\psi}{dt}  = \left[- \frac{\hbar^2}{2m}\frac{\partial^2}{\partial x^2}+V\right]\psi
\end{equation}
It is, then, not clear at all what the all the differentials actually are in a framework that implies countabilism. So, the present situation is that if we accept a mathematical foundation that implies countabilism, we cannot do physics. Consequently, from the applied side a mathematical foundation implying countabilism is not acceptable until all these problems have been solved---and they may be unsolvable.}

{From this point of view, the present framework of $\mathfrak{T}$ is preferable over any mathematical framework that implies countabilism. The reason is obvious: the framework of $\mathfrak{T}$ allows the uncountable set of the real numbers. The present framework of $\mathfrak{T}$ can therefor serve as a foundational framework for mathematics that is applicable to the sciences.}

\subsection{Comparison to foundational approaches based on type theory}\label{sect:TypeTheory}
{Over the years approaches to a foundation for mathematics based on Russel's type theory have been developed, such as
\begin{enumerate}[(i)]
  \item New Foundations (NF), a finitely axiomatizable theory developed by Quine \cite{Quine};
  \item intuitionistic type theory (MLTT) developed by Martin-L\"{o}f \cite{Martin,Martin2}, which led to the development of Homotopy Type Theory (HoTT) \cite{Awodey};
  \item NF with urelements (NFU), a finitely axiomatized variant of NF developed by Jensen \cite{Jensen}
\end{enumerate}
For a discussion on foundational matters involving these approaches, see \cite{ARDM}. The section ``Is the mathematical brain self-stratifying?'' therein shows that results from the Belgian school of NF-istes explain why mathematicians solely interested in geometry will think in stratified terms and find no use for abstract recursion; and therefore find versions of Zermelo set theory (Z) or the even weaker MacLane set theory sufficient for the mathematics they want to do. It has been shown, however, that Z does not support the definition of forcing, while Z extended with provident set theory does and is equiconsistent with Z \cite{ARDM2}. For a substantial discussion of MacLane set theory, see \cite{ARDM3,MacLane,ARDM4}.}

{Of the above approaches, NF has been shown to be finitely axiomatizable \cite{Hailperin}, and its consistency seems to have been proven recently \cite{Wilshaw}. But it has been proved that NF is inconsistent with AC \cite{Specker}. Consequently, to quote Holmes, ``Quine's `New Foundations has a bad reputation ... No one has been able to derive a contradiction from NF, but the failure of AC makes it an unfriendly theory to work in'' \cite{Holmes}. In the framework of $\mathfrak{T}$, on the other hand, AC can be derived from the axioms of $\mathfrak{T}$ \cite{CabboletAC}. {Now AC is indispensable for foundational work in mathematics: it not only allows to prove the existence of infinite sets which cannot be explicitly constructed, such as a basis for infinite-dimensional vector spaces and a well-ordering of the reals, but it is also essential for the proofs of basic theorems in various mathematical disciplines, such as Tychonoff's theorem in topology stating that the product of compact spaces is compact.} Therefore, $\mathfrak{T}$ is preferable over NF as a foundational theory for mathematics.}

{As to MLTT and HoTT, these are theoretical constructs that certainly have their strengths, but their intrinsic complexity undermines their role as a foundation for mathematics. In particular, it has been argued that the complexity of HoTT---with features like higher inductive types, univalence, and universe hierarchies---makes it less practical as a foundation for mainstream mathematics. For example, D\u{z}amonja argued that HoTT introduces unnecessary technical and conceptual overhead without any clear foundational advantages over ZFC \cite{Dzamonja}. Furthermore, working in the framework of HoTT requires a solid understanding of both type theory and homotopy theory: this can be a practical barrier for the use of HoTT as a foundation for mathematics in the sciences. In particular, it seems unthinkable that engineers will ever use the framework of HoTT to formalize differential equations that come up in the mathematical modelling of real world problems: this is a complex undertaking compared to more established foundational frameworks. In addition, Ladyman and Presnell argued that HoTT is not self-contained, as it leans on external concepts from homotopy theory \cite{Ladyman}. So, $\mathfrak{T}$ is preferable over MLTT and HoTT as a foundational theory for mathematics, because it is a self-contained theory consisting of simple axioms for which there are no practical barriers with regard to their use in the sciences.}

{Last but not least, NFU contains a nonconstructive axiom, its Axioms of Atoms, by which the existence of urelements (atoms) is assumed. These urelements are further undefined: the philosophy behind that is that these urelements can be objects of the physical universe, so that physical objects hereby become available in a mathematical discussion \cite{Holmes}. We reject this idea as preposterous: the purpose of mathematics is to provide an exact language for the sciences, but this language is in the first place a \emph{written} language. So whenever we write down a singleton of a physical object---e.g. a chair---as $\{\phi\}$, then in any case the symbol `$\phi$' is $^*$\textbf{not}$^*$ the physical object \emph{because the physical object cannot be written down}. At best the written symbol $\phi$ in this singleton $\{\phi\}$ is a mathematical thing that \emph{uniquely refers} to a physical object. But such constructions are possible without assuming the existence of urelements that can be physical objects. A second argument against NFU is that it is less user-friendly: if we have a set $X$ and we want to construct a subset $Y\subset X$ by applying separation, as in $Y = \{z\in Z\ |\ \Phi(z)\}$, then the existence of $Y$ is only guaranteed if $\Phi$ is a stratified formula---this is something that needs to be checked before separation can be applied. In the framework of $\mathfrak{T}$ on the other hand, we can apply separation as in the framework of ZFC. For the above reasons, $\mathfrak{T}$ is preferable over NFU as a foundational theory for mathematics.}

\subsection{Conclusion}\label{sect:conclusion}
Summarizing, we have proved that $\mathfrak{T}$ is relatively consistent with ZF, as in $\left[ \nvdash_{ZF} \bot\right] \Rightarrow \left[  \not\Vvdash_{\mathfrak{T}} \bot\right]$. Note that this is a relation between two theories in different logical frameworks! Given that $\mathfrak{T}$ has been proved to be stronger than ZF, we therefore tentatively conclude that we have made a concrete step towards proving that $\mathfrak{T}$ is an advancement in the foundations of mathematics.

{In addition, we have given arguments as to why $\mathfrak{T}$ is preferable over other foundational theories for mathematics. In particular, we prefer $\mathfrak{T}$ over any foundational theory reflecting the adage ``set theory is obsolete'', put forward by MacLane as mentioned in \cite{ARDM3,MacLane}: $\mathfrak{T}$ is a marriage of set theory and category theory, and we are not planning a divorce. As an objection to $\mathfrak{T}$} one may say that it is doubtful that nonclassical axioms such as the sum function axiom SUM-F will ever be used in everyday mathematical practice, let alone in everyday practice in the mathematical sciences. The point, however, is that one does not have to use SUM-F. Since ZF can be derived from $\mathfrak{T}$, standard first-order theorems of $\mathfrak{T}$ suffice for everyday practice in mathematics and the mathematical sciences: the nonclassical part of $\mathfrak{T}$ stays in the background. The essential point is that if we begin mathematical work within the framework of $\mathfrak{T}$, we do not need to worry about the existence of a countable model, as it does not exist. {So, why not switch to $\mathfrak{T}$?}
%\ \\
%\newpage
\appendix
%\noindent {\bf \large APPENDIX}

\section*{Appendix: outline of Set Matrix Theory}\label{app:SMT}

Informally, SMT is a generalisation of ZF: in the universe of SMT not only sets exist, but also set matrices that are constructed from a finite number of sets---these set matrices are \emph{urelements} in the sense that they are not sets but can be elements of sets. The language $L_{SMT}$ of SMT is the language $L_{ZF}$ of ZF, with the constant $\emptyset$ and with (labeled) Roman variables $x, y, \ldots$ varying over sets, extended with
\begin{enumerate}[(i)]
  \item (labeled) Greek variables $\alpha,\beta, \ldots$ varying over set matrices;
  \item for each pair of positive integers $m$ and $n$, an $mn$-ary function symbols $f_{m\times n}$, with matrix notation
  $\left [ \begin{array}{ccc} t_{11} & \ldots & t_{1n} \\ \vdots &   & \vdots \\ t_{m1} & \ldots & t_{mn} \end{array} \right ]$ for $f_{m\times n}(t_{11}, \ldots, t_{mn})$.
\end{enumerate}
A constructive \textbf{set matrix axiom schema} states that given any $m\cdot n$ set matrices $\alpha_{11}, \ldots, \alpha_{mn}$, there exists an $m\times n$ set matrix $\beta$ whose entries are the $\alpha_{ij}$'s:
\begin{equation}\label{eq:SetMatrixAxiom}
\forall \alpha_{11} \cdots \forall \alpha_{mn} \exists \beta: \beta = \left [ \begin{array}{ccc} \alpha_{11} & \ldots & \alpha_{1n} \\ \vdots &   & \vdots \\ \alpha_{m1} & \ldots & \alpha_{mn} \end{array} \right ]
\end{equation}
The \textbf{reduction axiom} of SMT identifies a $1\times 1$ set matrix with its entry:
\begin{equation}\label{eq:reduction}
  \forall x : [x] = x
\end{equation}
This way, quantification over all set matrices implies quantification over all sets.

Likewise, the \textbf{omission axiom schema} identifies a $1\times 1$ set matrix whose entry is an $m\times n$ set matrix with its entry:
\begin{equation}\label{eq:omission}
\forall \alpha_{11} \cdots \forall \alpha_{mn} : \left [ \left [ \begin{array}{ccc} \alpha_{11} & \ldots & \alpha_{1n} \\ \vdots &   & \vdots \\ \alpha_{m1} & \ldots & \alpha_{mn} \end{array} \right ] \right ]= \left [ \begin{array}{ccc} \alpha_{11} & \ldots & \alpha_{1n} \\ \vdots &   & \vdots \\ \alpha_{m1} & \ldots & \alpha_{mn} \end{array} \right ]
\end{equation}
The \textbf{division axiom schema} of SMT then states that set matrices of different sizes are not identical; this schema includes axioms stating that a set matrix of any size other than $1\times 1$ is not a set:
\begin{equation}\label{eq:DivisionAxiom}
\forall x \forall \alpha_{11} \cdots \forall \alpha_{mn}: x \neq \left [ \begin{array}{ccc} \alpha_{11} & \ldots & \alpha_{1n} \\ \vdots &   & \vdots \\ \alpha_{m1} & \ldots & \alpha_{mn} \end{array} \right ] \ \ \ m\cdot n \geq 2
\end{equation}
And the \textbf{epsilon axiom schema} states that an $m\times n$ set matrix with $m\cdot n \geq 2$ has no elements in the sense of the $\in$-relation:
\begin{equation}\label{eq:EpsilonAxiom}
\forall \alpha_{11} \cdots \forall \alpha_{mn} \forall \beta: \beta \not\in \left [ \begin{array}{ccc} \alpha_{11} & \ldots & \alpha_{1n} \\ \vdots &   & \vdots \\ \alpha_{m1} & \ldots & \alpha_{mn} \end{array} \right ] \ \ m\cdot n \geq 2
\end{equation}
The \textbf{extensionality axiom schema for set matrices} then states that two $m\times n$ set matrices are equal if their entries are equal:
\begin{equation}
\forall \alpha_{11} \cdots \forall \alpha_{mn} \forall \beta_{11} \cdots \forall \beta_{mn}  : \left [ \begin{array}{ccc} \alpha_{11} & \ldots & \alpha_{1n} \\ \vdots &   & \vdots \\ \alpha_{m1} & \ldots & \alpha_{mn} \end{array} \right ] = \left [ \begin{array}{ccc} \beta_{11} & \ldots & \beta_{1n} \\ \vdots &   & \vdots \\ \beta_{m1} & \ldots & \beta_{mn} \end{array} \right ] \Leftrightarrow \alpha_{11} = \beta_{11} \wedge \ldots \wedge \alpha_{mn} = \beta_{mn}
\end{equation}
As to the set-theoretical axioms of SMT, these are generalizations of the axioms of ZF, since set matrices that are not sets can still be elements of sets. E.g. the \textbf{generalized extensionality axiom for sets} becomes
\begin{equation}\label{eq:GeneralizedEXT}
  \forall x \forall y : x = y \Leftrightarrow \forall \alpha (\alpha \in x \Leftrightarrow \alpha\in y)
\end{equation}
Without further explanation, the generalizations of the constructive axioms of ZF are the following:
\begin{gather}
\exists x : x = \emptyset \wedge \forall \alpha (\alpha \not \in x)\label{eq:GeneralizedEMPTY} \\
\forall x \exists y \forall \alpha : \alpha \in y \Leftrightarrow \alpha \in x \wedge \Phi(\alpha)\label{eq:SEP}  \\
\forall \alpha \forall \beta \exists x \forall \gamma : \gamma \in x \Leftrightarrow \gamma = \alpha \vee \gamma = \beta \\
\forall x : \forall \alpha (\alpha \in x \Rightarrow \exists u (u = \alpha )) \Rightarrow \exists y \forall \beta (\beta \in y \Leftrightarrow \exists z (z \in x \wedge \beta \in z)) \\
\forall x \exists y \forall \alpha : \alpha \in y \Leftrightarrow \exists u (u = \alpha \wedge u \subseteq x) \\
\exists x : \emptyset \in x \wedge \forall y (y \in x \Rightarrow \{y \} \in x) \\
\forall x : \forall \alpha\in x \exists!\beta \Phi(\alpha, \beta) \Rightarrow \exists y \forall \zeta (\zeta\in y \Leftrightarrow \exists \gamma : \gamma \in x \wedge \Phi(\gamma, \zeta))\label{eq:GeneralizedREP}
\end{gather}
In addition, SMT has the \textbf{set of set matrices axiom schema}:
\begin{equation}\label{eq:SetOfSetMatricesAxiom}
  \forall x \exists y \forall \alpha : \alpha \in y \Leftrightarrow
  \exists \beta_{11} \ldots \exists \beta_{mn}(\alpha = \left [ \begin{array}{ccc} \beta_{11} & \ldots & \beta_{1n} \\ \vdots &   & \vdots \\ \beta_{m1} & \ldots & \beta_{mn} \end{array} \right ] \wedge
\beta_{11} \in x \wedge \ldots \wedge \beta_{mn} \in x)
\end{equation}
These are all the constructive axioms of SMT.

\section*{Author contributions}
A.R.D. Mathias suggested the proof method. M.J.T.F. Cabbolet developed the proof. A.R.D. Mathias and M.J.T.F. Cabbolet wrote the paper.

\section*{Acknowledgement}
The authors would like to thank Harrie de Swart (Tilburg University) for his useful comments.

\end{document}